\documentclass[11pt]{article}
\usepackage{fullpage,graphicx,psfrag,latexsym,color}
\usepackage{algorithm,algorithmic}
\usepackage{dcolumn}
\newcolumntype{d}[1]{D{.}{.}{#1.#1}}
\usepackage{amsmath,tikz}
\usepackage{subfig}
\usetikzlibrary{shapes.multipart}
\newcommand{\BEAS}{\begin{eqnarray*}}
\newcommand{\EEAS}{\end{eqnarray*}}
\newcommand{\BEA}{\begin{eqnarray}}
\newcommand{\EEA}{\end{eqnarray}}
\newcommand{\BEQ}{\begin{equation}}
\newcommand{\EEQ}{\end{equation}}
\newcommand{\BIT}{\begin{itemize}}
\newcommand{\EIT}{\end{itemize}}

\newcommand{\ie}{{\it i.e.}}

\newcommand{\reals}{{\mbox{\bf R}}}
\newcommand{\symm}{{\mbox{\bf S}}}

\newcommand{\Tr}{\mathop{\bf tr}}

\renewcommand{\symm}{\mathbf{S}}  

\newcommand{\proj}{\mathcal{P}}
\newcommand{\ch}{\mathrm{ch}}
\newcommand{\new}[1]{N_{#1}}
\newcommand{\anc}[1]{A_{#1}}
\newcommand{\I}[1]{I_{#1}}
\newcommand{\Ip}[1]{J_{#1}}

\newcommand{\pa}[1]{\mathrm{par}(#1)}

\newcommand{\SV}{\symm^n_V} 
\newcommand{\SVp}{\symm^n_{V,+}}
 
\newcommand{\SVc}{\symm^n_{V,\mathrm{c}}}

\newcommand{\sdot}{\mathbin{\bullet}}

\newtheorem{theorem}{Theorem}

\newcounter{algorithmctr}[section]
\renewcommand{\thealgorithmctr}{\thesection.\arabic{algorithmctr}}
\newenvironment{algdesc}[1]    {\refstepcounter{algorithmctr}\begin{list}{}{        \setlength{\rightmargin}{.05\linewidth}        \setlength{\leftmargin}{.05\linewidth}}        \item[]{\setlength{\parskip}{0ex}\bigskip\par         \nopagebreak         \underline{{\bf Algorithm \thealgorithmctr.} \emph{#1.}}}}    {{\setlength{\parskip}{-1ex}\nopagebreak\smallskip\par} \end{list}}

\title{Logarithmic barriers for sparse matrix cones}
\author{Martin S.\ Andersen\thanks{
    Electrical Engineering Department, University of California, 
    Los Angeles. Email: martin.andersen@ucla.edu, vandenbe@ee.ucla.edu.} 
\and Joachim Dahl\thanks{MOSEK ApS, Fruebjergvej 3, 2100 K{\o}benhavn {\O}.
   Email: dahl.joachim@gmail.com.} 
\and Lieven Vandenberghe\footnotemark[1]}
\date{}

\begin{document}
\maketitle

\begin{abstract}
Algorithms are presented for evaluating gradients and Hessians of 
logarithmic barrier functions for two types of convex cones:
the cone of positive semidefinite matrices with a given sparsity pattern, 
and its dual cone, the cone of sparse matrices with the same pattern 
that have a positive semidefinite completion.
Efficient large-scale algorithms for evaluating these barriers and 
their derivatives are important in interior-point methods for 
nonsymmetric conic formulations of sparse semidefinite programs.  
The algorithms are based on the multifrontal method for sparse 
Cholesky factorization. 
\end{abstract}

\section{Introduction}
\subsection{Log-det barrier for sparse matrices}
We discuss algorithms for evaluating the gradient and 
Hessian of the `log-det' barrier 
\[
f(X) = -\log\det X
\]
when $X$ is large, sparse, and positive definite.
We take $f$ as a function from $\SV$ to $\reals$, where 
$V$ is the filled sparsity pattern of the Cholesky factor of $X$, and
$\SV$ denotes the set of symmetric matrices with sparsity pattern $V$.
With this convention, and for the standard trace inner product
$A\sdot B = \Tr(AB)$ of symmetric matrices, the gradient of $f$ at $X$ is
\BEQ \label{e-grad-def}
 \nabla f(X) = -\proj(X^{-1})
\EEQ
where $\proj$ denotes projection on $V$, \ie, $\proj(A)_{ij} = A_{ij}$ if 
the pattern $V$ has a nonzero in position $i,j$ and $\proj(A)_{ij} = 0$ 
otherwise.
The algorithms presented in this paper exploit properties of
filled sparsity patterns (which are also known as \emph{chordal} or 
\emph{triangulated} patterns) to compute the gradient directly from the 
Cholesky factor of $X$, without calculating the rest of the inverse. 

The Hessian of $f$, interpreted as a function from $\SV$ 
to $\reals$, is defined by
\BEQ
 \nabla^2f(X)[Y] = \left. \frac{d}{dt} \nabla f(X+tY) \right|_{t=0}
 = \proj(X^{-1}YX^{-1}) \quad \forall Y \in \SV.
\EEQ
We are interested in efficient methods for evaluating this expression,
possibly for multiple matrices $Y$ simultaneously, without computing 
the entire inverse of $X$ or the products $X^{-1}YX^{-1}$.
We also discuss methods for evaluating the inverse Hessian 
$\nabla^2 f(X)^{-1}[Y]$.

The function $f$ has an important role as a logarithmic barrier
function for the convex cone
\[
  \SVp = \{X\in\SV \mid X\succeq 0\},
\]
the cone of positive semidefinite matrices with sparsity pattern $V$. 

\subsection{Conjugate barrier}
Equally important is the corresponding dual barrier function
\BEQ \label{e-dual-barrier-def} 
 f_*(S) = \sup_{X\in\symm^n_V} \left(-S\sdot X - f(X)\right).
\EEQ
This is the conjugate or Legendre transform of $f$ applied to $-S$.
The function $f_*$ is a logarithmic barrier function for the dual cone 
of $\SVp$,  which contains the symmetric matrices with sparsity pattern 
$V$ that have a positive semidefinite completion:
\[
   \SVc = \proj(\symm^n_+) = \{ \proj(X) \mid X\succeq 0\}.
\]
The dual barrier $f_*(S)$ can be computed as 
$f_*(S) = \log\det \hat X - n$  where $\hat X$ is the maximizer in the 
definition of $f_*$, \ie, the solution of the nonlinear equation
$\nabla f(X) = -S$ or
\BEQ \label{e-Xhat-eq}
 \proj(X^{-1}) = S
\EEQ
with variable $X\in\SV$.
The gradient and Hessian of $f_*$ at $S$ also follow from the
maximizer $\hat X$ by applying standard properties of Legendre 
transforms:
\BEQ \label{e-legendre-rels}
 \nabla f_*(S) = -\hat X, 
 \qquad \nabla^2 f_*(S) = \nabla^2 f(\hat X)^{-1}.
\EEQ
For general sparsity patterns $V$, the maximizer $\hat X$ 
in~(\ref{e-dual-barrier-def}) needs to be computed by iterative methods.  
For filled patterns $V$, however, efficient direct algorithms exist.   
The algorithms discussed in this paper compute a Cholesky factorization 
of $\hat X$, given the matrix $S$, using a finite recursion that is very 
similar and comparable in cost to a Cholesky factorization.

There is an interesting connection between the dual barrier $f_*$
and the maximum determinant positive definite completion problem, which 
has been extensively studied in linear algebra \cite{GJSW:84,Lau:01}.
The optimization problem in~(\ref{e-dual-barrier-def}) is the Lagrange 
dual of the convex optimization problem
\BEQ \label{e-compl}
 \begin{array}{ll}
 \mbox{minimize} &  -\log\det Z - n \\
 \mbox{subject to} & \proj(Z) = S,
 \end{array}
\EEQ
with variable $Z\in\symm^n$. The primal and dual optimal solutions
$\hat X$ and $Z$ are related by the optimality condition $Z^{-1} = \hat X$.
The solution $\hat X$ of~(\ref{e-Xhat-eq}) is therefore the inverse 
of the maximum determinant positive definite completion of $S$. 

\subsection{Applications}
Efficient gradient and Hessian evaluations for $f$ and $f_*$ are 
critical to the performance of interior-point methods for conic 
optimization problems associated with the cones $\SVp$ and $\SVc$
\cite{ADV:10,SrV:04}.  
Consider the pair of primal and dual cone linear programs
(LPs)
\BEQ \label{e-cone-LPs}
\begin{array}[t]{ll}
\mbox{minimize} & c^T x  \\
\mbox{subject to} & \sum\limits_{i=1}^m x_i A_i + X = B \\
 & X \in \SVp
\end{array} \qquad\qquad
\begin{array}[t]{ll}
\mbox{maximize}  & -B \sdot S \\
\mbox{subject to} & A_i \sdot S + c_i = 0, \quad i=1,\ldots,m \\
                  & S \in\SVc
\end{array} 
\EEQ
with variables $x\in\reals^m$, $X, S\in\SV$, and problem parameters 
$A_i, B \in\SV$, $c\in\reals^m$.  
Cone LPs of this type have been studied in sparse semidefinite 
programming with the goal of exploiting aggregate sparsity in the 
coefficient matrices $A_i$ and $B$ \cite{SrV:04,Bur:03,ADV:10}.
Matrix completion techniques and chordal sparse matrix properties
were first applied to semidefinite programming algorithms by 
Fukuda \emph{et al.} \cite{FKMN:00} in a sparse implementation
of the HRVW/KSH/M primal-dual algorithm.

A primal barrier method for~(\ref{e-cone-LPs}) requires at each iteration
the evaluation of the gradient $\nabla f(X)$ and the solution
of a positive definite equation $H\Delta y = g$ with coefficients 
$H_{ij} = A_i \sdot (\nabla^2f(X)[A_j])$.   
Efficient techniques for evaluating the Hessian $\nabla^2f(X)[A_j]$ are 
therefore important in large-scale implementations.
A dual barrier method for the cone programs involves evaluations of 
the gradient $\nabla f_*(S)$ and a set of linear equations 
$H\Delta y = g$ with coefficients 
$H_{ij} =  A_i \sdot (\nabla^2 f_*(S)^{-1}[A_j])$.
From the relations~(\ref{e-legendre-rels}), we see that this 
requires the inverse $\hat X$ of the maximum determinant positive 
definite matrix completion of $S$ and the evaluation of the Hessian 
$\nabla^2 f(\hat X)[A_j]$ at $\hat X$.  
We refer the reader to~\cite{ADV:10,SrV:04} for more details.

The problem of computing a projected inverse $\proj(X^{-1})$ (and, more 
generally, computing a subset of the entries of the inverse of 
a sparse positive definite matrix) has also been studied in statistics 
\cite{GoP:80,ADRRU:10}.
Efficient algorithms for computing the gradient of $f$ are 
important in maximum likelihood estimation problems involving
Gaussian distributions, for example, in sparse inverse covariance 
selection \cite{DVR:08}.
Consider, for example, the covariance selection problem with 
$1$-norm penalty 
\[
\begin{array}{ll}
\mbox{minimize} & C\sdot X - \log\det X + \rho \|X\|_1,
\end{array}
\]
which has been studied by several authors
\cite{HLP+:06,BED:08,FHT:08,ABE:08,SMG:10,LiT:10}.
In this problem, $X$ is the inverse covariance matrix of a Gaussian
random variable and $C$ is a sample covariance.  The first two terms
in the objective form the negative log-likelihood function of $X$ (up
to constants), and the penalty term is added to promote sparsity 
in the solution $X$.  In problems of high dimension, it may
be unrealistic and impractical to regard $X$ as a dense matrix
variable.  Instead, one can start with a partially specified pattern
based on prior knowledge, and use the penalized covariance selection
to identify additional zeros.  The problem can then be posed as an
optimization problem over $\SV$, where $V$ is the known sparsity
pattern.  This greatly simplifies the cost of calculating the gradient
of the smooth terms in the objective and makes it possible to solve
very large covariance selection problems using first-order methods
that require the gradient of $-\log\det X$ at each iteration.

Other applications include matrix approximation problems with 
sparse positive definite matrices as variables, for example,
computing sparse quasi-Newton updates \cite{Fle:95,Yam:07}.

\subsection{Related work and outline of the paper}
We refer to the algorithms in this paper as \emph{multifrontal} and 
\emph{supernodal} because of their resemblance to multifrontal and 
supernodal multifrontal algorithms for Cholesky factorization 
\cite{DuR:83,Liu:92}.  The multifrontal Cholesky factorization 
is reviewed in Section~\ref{s-chol} and a supernodal variant, formulated
in terms of clique trees, is described in Section~\ref{s-supernodal}.

In Section~\ref{s-grad}, we introduce multifrontal algorithms for
computing the gradients of $f$ and $f_*$.  Similar algorithms for 
evaluating $\nabla f$ are discussed in~\cite{CaD:95,ADRRU:10}.
The close connection between the problem of computing the gradient of
$f(X) = -\log\det X$ and a Cholesky factorization is easily understood
from the chain rule of differentiation.  A practical method for
computing $f(X)$ will calculate a sparse Cholesky factorization of
$X$, for example, $X=LDL^T$ with $L$ unit lower triangular and $D$
diagonal, and then evaluate $f(X) = -\sum_i \log D_{ii}$.  By applying
the chain rule to a sparse factorization algorithm, the gradient of
$f(X)$ can be evaluated at essentially the cost of the factorization
itself.  Moreover, the differentiation can be automated using reverse
automatic differentiation software \cite{GrW:08}.
Although the algorithm in Section~\ref{s-primal-grad} 
can be obtained from the chain rule,
we give a straightforward direct derivation.  This not only simplifies
the notation and description of the algorithm, it also helps reduce
the memory requirements, which can be high in a straightforward
application of reverse differentiation because of the large number of
intermediate auxiliary matrices (update and frontal matrices)
generated during the factorization.  

As mentioned earlier, evaluating the gradient $\nabla f_*(S)$ is 
equivalent to inverting the mapping $X
\mapsto \nabla f(X)$, and an algorithm for evaluating $\nabla f_*$ is
therefore easily derived from the algorithm for $\nabla f$
(see Section~\ref{s-dual-gradient}).

In Section~\ref{s-hess}, we examine the problem of computing
the Hessians and inverse Hessians of $f$ and $f_*$, and more 
specifically,
the problem of evaluating expressions of the form $\nabla^2 f(X)[U]$
and $\nabla^2 f_*(S)[V]$.
Again, the algorithms follow conceptually from the chain rule  and 
can be obtained by applying automatic differentiation techniques.
An explicit description allows us to optimize the efficiency and 
memory requirements.  This is particularly important in the case 
of sparse arguments $U$, $V$.  As an important by-product, we 
define a factorization 
$\nabla^2 f(X) = \mathcal R^\mathrm{adj} \circ \mathcal R$
and present efficient methods for evaluating the factors
$\mathcal R$ and $\mathcal R^\mathrm{adj}$ separately.

The methods presented in the paper are closely related to
the barrier evaluation algorithms of~\cite{DVR:08,DaV:09a}.
These algorithms were formulated as recursions over clique trees 
and can be interpreted as supernodal versions of the multifrontal
algorithms presented in Sections~\ref{s-chol}--\ref{s-hess}.
We elaborate on the connections in Sections~\ref{s-cliques} 
and~\ref{s-supernodal}.
In contrast to the clique tree methods in~\cite{DVR:08,DaV:09a}, the 
algorithms described here work with data structures that are 
widely used in sparse Cholesky factorization algorithms  
(namely, the compressed 
column storage format and elimination trees).  As a result they are 
simpler to implement and more readily combined with techniques 
from the recent literature on sparse matrix factorization algorithms.
Some possible further improvements in this direction are 
mentioned in the conclusions (Section~\ref{s-conclusions}).

\subsection{Notation}   \label{s-notation}
We identify a symmetric sparsity pattern $V$ with the positions of its 
lower-triangular nonzeros.
In other words, a symmetric matrix $X$ has sparsity pattern $V$ if 
$X_{ij} = X_{ji} = 0$ for $(i,j)\not\in V$.
The entries $X_{ij}$ and $X_{ji}$ for $(i,j)\in V$ are treated as 
(structurally) nonzero, although they are allowed to be numerically zero.
A lower-triangular martrix $L$ has sparsity pattern $V$ 
if the symmetric matrix $L+L^T$ has sparsity pattern $V$.

An \emph{index set} is a sorted subset of the 
integers $\{1,2,\ldots,n\}$.  
The number of elements in the index set $I$ is denoted $|I|$ and its 
$k$th element $I(k)$.
If $I$ and $J$ are two index sets with $J\subset I$, we define an 
$|I| \times |J|$ matrix $E_{IJ}$ with entries
\[
  (E_{IJ})_{ij} = \left\{ \begin{array}{ll}
    1 & I(i) = J(j)  \\ 0 & \mbox{otherwise.} \end{array}\right.
\]
This notation will be used in expressions $E_{IJ}^T BE_{IJ}$ and
$E_{IJ}CE_{IJ}^T$, which have the following meaning:
if $A$ is a symmetric matrix of order $n$ and $B=A_{II}$ is the 
principal submatrix indexed by $I$, then the matrix $E_{IJ}^TBE_{IJ}$ is 
equal to $A_{JJ}$, the principal submatrix indexed by $J$.
The adjoint operation $B = E_{IJ} C E_{IJ}^T$, applied to a symmetric 
matrix $C$ of order $|J|$, can be interpreted as first embedding $C$ 
as the $|J|\times |J|$-block $A_{JJ}$ of an otherwise zero 
$n\times n$ matrix $A$, and then extracting the $|I|\times |I|$
submatrix $B = A_{II}$.

\section{Elimination trees} \label{s-etree}
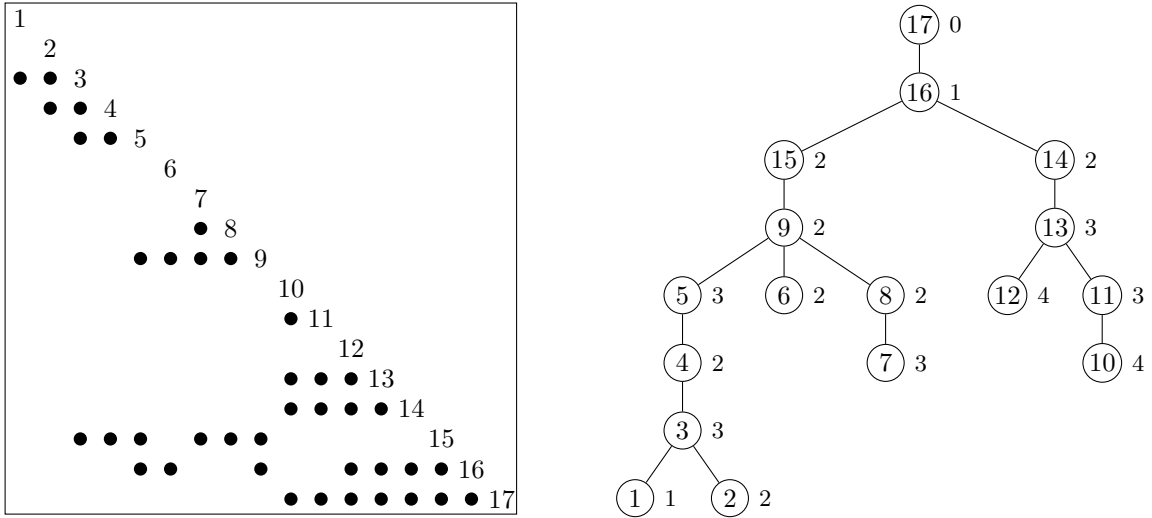
\begin{figure}

\hspace*{\fill}
\begin{tikzpicture}[scale=0.400000]
   \tikzset{VertexStyleA/.style = {shape = circle, minimum size = 5pt, 
       inner sep = 0 pt, fill=black}}

   \foreach \i/\j in {
       2/0,
       2/1, 3/1,
       3/2, 4/2, 14/2,
       4/3, 14/3,
       8/4, 14/4, 15/4, 
       8/5, 15/5,
       7/6, 8/6, 14/6,  
       8/7, 14/7, 
       14/8, 15/8,
       10/9, 12/9, 13/9, 16/9,
       12/10, 13/10, 16/10,
       12/11, 13/11, 15/11, 16/11,
       13/12, 15/12, 16/12,
       15/13, 16/13,
       15/14, 16/14,
       16/15 
      }
      \node[VertexStyleA] at (\j,-\i){};

   \foreach \k/\i in {
       1/0, 2/1, 3/2, 4/3, 5/4, 6/5, 7/6, 8/7, 9/8, 10/9,
       11/10, 12/11, 13/12, 14/13, 15/14, 16/15, 17/16} 
      \node[font=\small] at (\i,-\i) {$\k$};

   \draw (-0.5,0.5) rectangle (16.5,-16.5);
\end{tikzpicture}\hspace*{\fill}
\hspace*{\fill}
\begin{tikzpicture}[scale=9.00]
      \tikzset{VertexStyle/.style = {shape = circle, draw,  minimum size = 14, inner sep = 1, fill=none}}

      \foreach \l/\k/\d/\x/\y in {
       1/1/1/0.13/0.0,
       2/2/2/0.27/0.0,
       3/3/3/0.2/0.1,
       4/4/2/0.2/0.2,
       5/5/3/0.2/0.3,
       6/6/2/0.35/0.3,
       7/7/3/0.5/0.2,
       8/8/2/0.5/0.3,
       9/9/2/0.35/0.4,
       10/10/4/0.82/0.2,
       11/11/3/0.82/0.3,
       12/12/4/0.68/0.3,
       13/13/3/0.75/0.4,
       14/14/2/0.75/0.5,
       15/15/2/0.35/0.5,
       16/16/1/0.55/0.6,
       17/17/0/0.55/0.7}
       \node[VertexStyle, font=\small](\k) 
           [label = right: \footnotesize $\d$] at (\x,\y){\l};

      \foreach \i/\j in {
       1/3,
       2/3,
       3/4,
       4/5,
       5/9,
       6/9,
       7/8,
       8/9,
       9/15,
       15/16,
       10/11,
       11/13,
       12/13,
       13/14,
       14/16,
       16/17}
       \draw (\i)--(\j);
\end{tikzpicture}
\hspace*{\fill}
\caption{Filled pattern and elimination tree.  The numbers next to 
 vertices in the elimination tree are the monotone degrees $|\I{k}|$. 
\label{f-etree1}}
\end{figure}

This section provides some background on sparse matrices and elimination 
trees \cite{Liu:90,Dav:06}.
We define a Cholesky factorization as a factorization 
\[
 X = LDL^T
\]
with $L$ unit lower-triangular and $D$ positive diagonal.
The sparsity pattern $V$ of the Cholesky factor $L$ has the following 
fundamental property:  
\BEQ \label{e-filled}
  i > j > k, \quad (i,k) \in V, \quad (j,k) \in V\quad 
\Longrightarrow \quad (i,j)\in V.
\EEQ
(In other words, excluding accidental cancellation, 
$L_{ik} \neq 0$ and $L_{jk} \neq 0$ implies $L_{ij} \neq 0$.)
This property distinguishes a Cholesky factor from a general sparse 
lower-triangular matrix.
In the example in Figure~\ref{f-etree1}, the presence of nonzeros in 
positions $(5,3)$ and $(15,3)$ implies that the entry $(15,5)$ is nonzero.
The nonzeros in positions $(9,5)$, $(15,5)$, and $(16,5)$ imply that the 
entries in positions $(15,9)$, $(16,9)$, and $(16,15)$ are nonzero.

We use the notation $\I{k}$ to denote the sorted set of row indices 
of the nonzero entries below the diagonal in column $k$ of $L$.  
We also define $\Ip{k} = \I{k} \cup \{k\}$.  
The property~(\ref{e-filled}) implies that $\I{k}$ defines a complete 
subgraph of the filled graph, \ie, the matrix $L_{\Ip{k} \Ip{k}}$ is a 
dense lower-triangular matrix.
For example, it can be verified that the submatrix
indexed by $\Ip{5} = \{5,9,15,16\}$ in Figure~\ref{f-etree1} is dense. 
The number of nonzeros below the diagonal in column $k$, \ie, the 
cardinality $|\I{k}|$ of $\I{k}$, is called the \emph{monotone degree}  of 
vertex  $k$.

The \emph{elimination tree} (etree) is defined in terms of the sparsity 
pattern of the factor $L$ as follows.  
It is a tree (or a forest if $L$ is reducible) with $n$ vertices, 
labeled $1$ to $n$.
The parent of vertex $k$ is the row index $j$ of the first nonzero
below the diagonal of column $k$ of $L$, \ie, the vertex $\min \I{k}$.
As a consequence, each vertex has a lower index than its parent,
so the vertices in the elimination tree are numbered in a 
\emph{topological ordering}. 
An example is shown in Figure~\ref{f-etree1}.

We will use two important properties of elimination trees. 
\begin{theorem} \emph{\cite[theorem 3.1]{Liu:90}} \label{t-ancestor}
If $j \in \I{k}$, then $j$ is an ancestor of $k$ 
in the elimination tree.  
\end{theorem}
Note that the converse does not hold.  

\begin{theorem} 
\emph{\cite[theorem 3.1]{Liu:92}}  
\label{t-column-dep}
If vertex $j$ is an ancestor of vertex $k$ in the elimination tree, 
then the nonzero structure of $(L_{jk}, L_{j+1,k}, \ldots, L_{nk})$
is contained in the structure of $(L_{jj}, L_{j+1,j}, \ldots, L_{nj})$.
\end{theorem}
In the notation for the column structure introduced above, 
this theorem asserts that
if $j$ is an ancestor of $k$, then
\[
 \I{k} \cap \{j, j+1, \ldots, n\} \subseteq \Ip{j}.
\]
In particular, if $j$ is the parent of $k$ (hence, by definition, $j$ is 
the first element of $\I{k}$ and  therefore
$\I{k} \subseteq \{j,j+1,\ldots,n\}$), then
\BEQ \label{e-etree-ch}
 \I{k} \subseteq \Ip{j}, \qquad |\I{k}| \leq |\I{j}| + 1.
\EEQ
By applying these inequalities recursively to a path 
$i_1$, $i_2$, \ldots, $i_r$ from a vertex $i_1$ in the elimination tree
to one of its  ancestors $i_r$, we obtain a chain of inclusions
\BEQ \label{e-etree-path}
 \I{i_1} \; \subseteq \; \Ip{i_2} 
 \;  \subseteq \; \{i_2\} \cup \Ip{i_3}  \; \subseteq
 \;  \cdots \; \subseteq \; \{i_2, i_3, \ldots, i_{r-1} \} \cup \Ip{i_r} 
\EEQ
and inequalities
\BEQ \label{e-tree-path-I}
 |\I{i_1}| \; \leq \; |\I{i_2}| + 1 \; \leq \; 
 |\I{i_3}| + 2 \; \leq \; \cdots \; \leq \; |\I{i_r}| + r-1.
\EEQ
This can be verified in the elimination tree in Figure~\ref{f-etree1}.  
As we move along a path from a leaf vertex to 
the root of the tree, the monotone degrees can increase or decrease, 
but they never decrease by more than one per step.

\section{Cholesky factorization and multiplication} \label{s-chol}
In this section, we review the multifrontal algorithm for Cholesky 
factorization \cite{DuR:83,Liu:92}.
We then describe a similar algorithm for the related problem of 
computing a matrix, given its Cholesky factors.

\subsection{Cholesky factorization}
Recall that we define the Cholesky factorization as a decomposition 
$X = LDL^T$, with $D$ positive diagonal and $L$ unit lower-triangular.
The formulas for $L$ and $D$ are easily derived from the equation 
$X=LDL^T$.
The $\Ip{j} \times \Ip{j}$ block of the factorization is
\BEAS
\lefteqn{
 \left[\begin{array}{cc} 
 X_{jj} & X_{\I{j}j}^T \\ X_{\I{j}j} & X_{\I{j}\I{j}} 
 \end{array}\right] } \\
 &  = &  
 \sum_{k<j} D_{kk} \left[\begin{array}{c} 
 L_{jk} \\ L_{\I{j}k} \end{array} \right]
 \left[\begin{array}{c} L_{jk} \\ L_{\I{j}k} \end{array} \right]^T 
 +
 D_{jj} \left[\begin{array}{c} 1 \\ L_{\I{j}j} \end{array}\right]
   \left[\begin{array}{c} 1 \\ L_{\I{j}j} \end{array}\right]^T
 + \sum_{k>j}D_{kk}
   \left[\begin{array}{c} 0 \\ L_{\I{j}k} \end{array}\right]
   \left[\begin{array}{c} 0 \\ L_{\I{j}k} \end{array}\right]^T.
\EEAS
The first column of the equation is
\BEQ \label{e-chol-a}
 \left[\begin{array}{cc} X_{jj} \\ X_{\I{j}j} \end{array}\right]
 = \sum_{k<j} D_{kk} L_{jk} \left[\begin{array}{c} L_{jk} \\ L_{\I{j}k}
 \end{array}\right]
 + D_{jj} \left[\begin{array}{c} 1 \\ L_{\I{j}j} \end{array}\right].
\EEQ
The multifrontal algorithm takes advantage of properties of the elimination
tree associated with $L$ to compute the sum on the right-hand side.
First, we recall that $L_{jk} \neq 0$ only if $k$ is a descendant 
of $j$ in the elimination tree (Theorem~\ref{t-ancestor}).  
The sum in~(\ref{e-chol-a}) can therefore be replaced by a 
sum over the proper descendants of $j$.
The set of proper descendants of vertex $j$ is
\[
T_j \setminus \{j\} = \bigcup_{i\in\ch(j)} T_i,
\]
where $T_j$ is the subtree of the elimination tree rooted at
vertex $j$ and $\ch(j)$ are the children of vertex $j$.
The equation~(\ref{e-chol-a}) then becomes
\BEQ \label{e-chol}
 \left[\begin{array}{cc} X_{jj} \\ X_{\I{j}j} \end{array}\right]
 = \sum_{i\in\ch(j)} \sum_{k\in T_i} D_{kk} L_{jk} 
  \left[\begin{array}{c} L_{jk} \\ L_{\I{j}k} \end{array}\right]
 + D_{jj} \left[\begin{array}{c} 1 \\ L_{\I{j}j} \end{array}\right].
\EEQ
Second, suppose that for each vertex $i$ in the elimination tree, we 
define a dense matrix 
\BEQ \label{e-Ui}
   U_i = -\sum_{k\in T_i} D_{kk} L_{\I{i}k} L_{\I{i}k}^T.
\EEQ
The matrix $U_i$ is called the \emph{update matrix} for vertex $i$. 
Using 
the definition of $E_{IJ}$ in Section~\ref{s-notation}, we can write
\BEQ \label{e-Ui-j}
 -\sum_{k\in T_i} D_{kk} 
  \left[\begin{array}{c} L_{jk}\\ L_{\I{j}k}\end{array}\right]
  \left[\begin{array}{c} L_{jk}\\ L_{\I{j}k}\end{array}\right]^T
   = E_{\Ip{j}\I{i}} U_i E_{\Ip{j}\I{i}}^T.
\EEQ
This follows from Theorem~\ref{t-column-dep}: if $i\in\ch(j)$
and $k\in T_i$, then 
\[
 \I{k} \cap \{j,j+1,\ldots,n\} \subseteq \I{i}
 \subseteq \Ip{j}
\]
and therefore $L_{\Ip{j}k} = E_{\Ip{j}\I{i}} L_{\I{i}k}$.
(The multiplication with $E_{\Ip{j}\I{i}}$ copies the entries
$L_{\I{i}k}$ to the correct position in $L_{\Ip{j}k}$ and inserts 
zeros for the other entries.)
Adding the first columns of 
$E_{\Ip{j}\I{i}} U_i E_{\Ip{j}\I{i}}^T$ for all $i\in\ch(j)$
therefore gives the first term on the right-hand side of~(\ref{e-chol}).

Thus, by combining the nonzero lower-triangular entries in column $j$ of 
$X$ (the left-hand side of~(\ref{e-chol}))
and the update matrices of the children of vertex $j$
(to assemble the sum on the right-hand side), we collect all 
the information needed to compute $D_{jj}$ and $L_{\I{j}j}$
from~(\ref{e-chol}).

Furthermore, the same equation~(\ref{e-chol}) shows how the update matrix 
$U_j$ for vertex $j$ can be calculated.  This is clearer if we 
rewrite~(\ref{e-chol}) in matrix form using~(\ref{e-Ui-j}) as
\BEA
 \left[\begin{array}{cc} 
  X_{jj} & X_{\I{j}j}^T \\ X_{\I{j}j} & 0 \end{array}\right]
   + \sum_{i\in\ch(j)} E_{\Ip{j}\I{i}} U_i E_{\Ip{j}\I{i}}^T 
  & = &
   D_{jj} \left[\begin{array}{c} 1 \\ L_{\I{j}j} \end{array}\right]
  \left[\begin{array}{c} 1 \\ L_{\I{j}j} \end{array}\right]^T
 + \left[\begin{array}{cc} 0 & 0 \\ 0 & U_j \end{array}\right] 
\label{e-mf-chol-0} \\
 & = & \left[\begin{array}{cc} 1 & 0 \\ L_{\I{j}j} & I \end{array}\right]
 \left[\begin{array}{cc} D_{jj} & 0 \\ 0 & U_j  \end{array}\right]
 \left[\begin{array}{cc} 1 & L_{\I{j}j}^T \\ 0 & I \end{array}\right].
\label{e-mf-chol}
\EEA
The first column of this equation is identical to~(\ref{e-chol}).
The 2,2 block follows from the definition of $U_j$ and the 
identity~(\ref{e-Ui-j}).
The matrix on the left-hand side of~(\ref{e-mf-chol}) is called the 
$j$th \emph{frontal matrix}.
The equation~(\ref{e-mf-chol}) shows that once the  frontal matrix
has been assembled, we can compute $D_{jj}$, $L_{\I{j}j}$, 
and $U_j$ by a pivot step.  

The resulting algorithm to compute $L$, $D$, given a positive definite $X$,
is summarized below.

\begin{algdesc}{Cholesky factorization} \label{a-chol}
\begin{list}{}{}
\item[\textbf{Input.}] A positive definite matrix $X$.
\item[\textbf{Output.}] The factors $L$, $D$ in the Cholesky 
factorization $X = LDL^T$.
\item[\textbf{Algorithm.}] 
Iterate over $j\in \{1,\ldots,n\}$ in topological order
(\ie, visiting each vertex of the elimination tree before its parent).
For each $j$, form the frontal matrix
\BEQ \label{e-mf-chol-1}
 F_j = \left[\begin{array}{cc}
  F_{11} & F_{21}^T \\ F_{21} & F_{22} \end{array}\right] 
 = \left[\begin{array}{cc} 
 X_{jj} & X_{\I{j}j}^T \\ X_{\I{j}j} & 0 \end{array}\right]
 + \sum_{i \in \ch(j)} E_{\Ip{j}\I{i}} U_i E_{\Ip{j}\I{i}}^T
\EEQ
and calculate $D_{jj}$, the $j$th column of $L$, and the 
$j$th update matrix $U_j$ from 
\BEQ \label{e-mf-chol-2}
 D_{jj} = F_{11}, \qquad  
 L_{\I{j}j} = \frac{1}{D_{jj}} F_{21}, \qquad 
 U_j = F_{22} - D_{jj}L_{\I{j}j}L_{\I{j}j}^T.
\EEQ
\end{list}
\end{algdesc}

In a practical implementation, with the lower-triangular part of $X$ stored
in a sparse format (typically, the compressed column structure or CCS;
see \cite{Dav:06}), one can overwrite $X_{jj}$ with $D_{jj}$ 
and $X_{\I{j}j}$ with $L_{\I{j}j}$ after cycle $j$.
The auxiliary matrices $F_j$ and $U_i$ are stored as dense matrices
(either as two separate arrays or by letting $U_j$ overwrite
the $2,2$ block of $F_j$).  The main step in the algorithm
is the level-2 BLAS operation in the calculation of $U_j$
in~(\ref{e-mf-chol-2}) \cite{DDHH:88}.
The frontal matrix $F_j$ can be discarded after the vertex $j$
has been processed. The update matrix $U_j$ can be discarded after the 
parent of vertex $j$ has been processed.  

\subsection{Cholesky multiplication}
The equation~(\ref{e-mf-chol-0}) 
also shows how the $j$th column of $X$ can 
be computed from $D_{jj}$, column $j$ of $L$, and the update matrices for 
the children of vertex $j$.
This yields an algorithm for the inverse operation of the Cholesky 
factorization, \ie, the matrix multiplication $LDL^T$, 
which will be important in Section~\ref{s-grad}.

\newpage
\begin{algdesc}{Cholesky product}   \label{a-chol-prod}
\begin{list}{}{}
\item[\textbf{Input.}] Cholesky factors $L$, $D$.
\item[\textbf{Output.}] The matrix $X=LDL^T$.
\item[\textbf{Algorithm.}] Iterate over $j\in\{1,\ldots,n\}$ in 
topological order. 
For each $j$, calculate $X_{jj}$, $X_{\I{j}j}$, and $U_j$ from
\BEQ \label{e-mf-prod}
 \left[\begin{array}{cc} 
 X_{jj} & X_{\I{j}j}^T \\ X_{\I{j}j} & -U_j \end{array}\right]
 = D_{jj} \left[\begin{array}{c} 1 \\ L_{\I{j}j} \end{array}\right]
       \left[\begin{array}{c} 1 \\ L_{\I{j}j} \end{array}\right]^T
  - \sum_{i \in \ch(j)} E_{\Ip{j}\I{i}}U_iE_{\Ip{j}\I{i}}^T.
\EEQ
\end{list}
\end{algdesc}

The matrices $U_i$ can be deleted after the parent of vertex $i$ has been 
processed.  In a practical implementation, we compute the left-hand side 
of~(\ref{e-mf-prod}) as a dense matrix, via a level-2 BLAS operation
for the outer-product on the right-hand side.
Then $X_{jj}$ and $X_{\I{j}j}$ are copied to the CCS structure for $X$. 

\section{Gradients} \label{s-grad}
In this section we describe `multifrontal' algorithms for
evaluating the gradients of the barrier and the dual barrier.

Recall that the primal gradient is defined as 
$\nabla f(x) = -\proj(X^{-1})$, where $\proj$ denotes projection on 
the filled pattern $V$ of $X$.
It is straightforward to show that the entries of the projected 
inverse $\proj(X^{-1})$ can be computed directly from the Cholesky 
factors $L$, $D$ without calculating any entries of $X^{-1}$ outside $V$ 
(see Section~\ref{s-primal-grad}).
This observation is the basis of several algorithms published in 
the literature.
The algorithm we describe here is equivalent to the \emph{inverse 
multifrontal algorithm} in \cite{CaD:95}, but we give a different
and shorter derivation. 
The projected inverse algorithm in~\cite{DVR:08} can be viewed as 
a supernodal variant of the algorithm discussed here
(see the discussion in Section~\ref{s-supernodal}).
Another closely related algorithm is described by Amestoy \emph{et al.}\
\cite{ADRRU:10}, who consider the problem of computing a few entries 
of the inverse of a large sparse matrix.

Evaluation of the dual gradient corresponds to the inverse operation,
\ie, the problem of solving the nonlinear equation $\proj(X^{-1}) = S$
with variable $X\in\SV$.
A multifrontal algorithm for this problem is derived in 
Section~\ref{s-dual-gradient}.

\subsection{Primal gradient} \label{s-primal-grad}
Define $Z = X^{-1}$ and $S=\proj(Z)$.  We are interested in an 
efficient method for computing $S$ from the Cholesky factors $L$ and $D$
of $X$.  The matrix $Z = L^{-T}D^{-1}L^{-1}$ satisfies
\[
  Z L = L^{-T} D^{-1}.
\]
The $\Ip{j} \times j$ block of this equation only involves entries 
of $Z$ in the projection $S = \proj(Z)$:  
\BEQ \label{e-ldl-inv}
 \left[\begin{array}{cc} 
  S_{jj}  & S_{\I{j}j}^T \\ 
  S_{\I{j}j} & S_{\I{j}\I{j}}
 \end{array}\right]
 \left[\begin{array}{c} 1 \\ L_{\I{j}j} \end{array} \right] = 
 \left[\begin{array}{c} 1/D_{jj} \\ 0 \end{array}\right].
\EEQ
The vertices of $\I{j}$ are ancestors of vertex $j$
(Theorem~\ref{t-ancestor}).
Therefore, if we calculate the columns of $S$ following a
reverse topological order of the vertices of the elimination tree,
then the matrix $S_{\I{j}\I{j}}$ is known when we arrive at column $j$.
Given $S_{\I{j}\I{j}}$, it is easy to compute $S_{\I{j}j}$ and $S_{jj}$
from~(\ref{e-ldl-inv}):
\BEQ \label{e-grad}
 S_{\I{j}j} = - S_{\I{j}\I{j}} L_{\I{j}j}, \qquad 
 S_{jj} = \frac{1}{D_{jj}} - S_{\I{j}j}^TL_{\I{j}j}.
\EEQ
Accessing the vectors $L_{\I{j}j}$ and $S_{\I{j}j}$ is easy 
if $L$ and the lower-triangular part of $S$ are stored in a CCS data 
structure. 
Retrieving $S_{\I{j}\I{j}}$ from the CCS representation of $S$
can be avoided by using an idea similar to the multifrontal Cholesky 
algorithm.
For each vertex $j$ of the elimination tree, we define a dense 
`update matrix'
\[
    V_j = S_{\I{j}\I{j}}.
\]
It follows from the properties of the elimination tree 
(theorem~\ref{t-column-dep} and equation~(\ref{e-etree-ch}))
and the definition
of $E_{\Ip{j}\I{i}}$ that if $i$ is a child of vertex $j$, then
\[
  V_i = E_{\Ip{j}\I{i}}^T 
  \left[\begin{array}{cc} 
   S_{jj} & S_{\I{j}j}^T \\ S_{\I{j}j} & V_j \end{array}\right] 
   E_{\Ip{j}\I{i}}.
\]
By using this formula to propagate $V_i$, we obtain a `multifrontal'
algorithm for computing $\proj(X^{-1})$.

\begin{algdesc}{Projected inverse}   \label{a-projinv}
\begin{list}{}{}
\item[\textbf{Input.}] The Cholesky factors $L$, $D$ of a positive
definite matrix $X$.
\item[\textbf{Output.}] 
The projected inverse $S = \proj(X^{-1}) = -\nabla f(X)$.
\item[\textbf{Algorithm.}]
Iterate over $j=\{1,2,\ldots,n\}$ in reverse topological order
(\ie, visiting each vertex before its children).
For each $j$, calculate $S_{jj}$ and $S_{\I{j}j}$ from
\BEQ \label{e-mf-gradient-1}
 S_{\I{j}j} = -V_j L_{\I{j}j}, \qquad 
 S_{jj} = \frac{1}{D_{jj}} - S_{\I{j}j}^TL_{\I{j}j},
\EEQ
and compute the update matrices  
\BEQ \label{e-mf-gradient-2}
  V_i = E_{\Ip{j}\I{i}}^T \left[\begin{array}{cc} 
 S_{jj} & S_{\I{j}j}^T \\ S_{\I{j}j} & V_j \end{array}\right] 
  E_{\Ip{j}\I{i}}, \quad
 i\in \ch(j). 
\EEQ
\end{list}
\end{algdesc}
The matrix $V_j$ can be discarded after cycle $j$, and $S_{\I{j}j}$ and 
$S_{jj}$ can overwrite $L_{\I{j}j}$ and $D_{jj}$ in a CCS data structure.
As in Algorithm~\ref{a-chol}, the algorithm involves operations
with the dense matrices $V_j$ and 
\[
  \left[\begin{array}{cc} S_{jj} & S_{\I{j}j}^T \\ 
  S_{\I{j}j} & S_{\I{j}\I{j}} \end{array}\right].
\] 
The main calculation is a level-2 BLAS operation (the matrix-vector 
product $V_jL_{\I{j}j}$).

\subsection{Dual gradient} \label{s-dual-gradient}
As mentioned in the introduction, the solution $X\in\SV$ of 
the equation $\proj(X^{-1}) = S$ is the inverse
of the maximum determinant positive definite completion of $S$,
and it is also the negative of $\nabla f_*(S)$.
The Cholesky factors of $X$ can be computed by solving for 
$L_{I_jj}$ and $D_{jj}$ from~(\ref{e-grad}) as in the following algorithm.

\begin{algdesc}{Matrix completion}   \label{a-completion}
\begin{list}{}{}
\item[\textbf{Input.}] A matrix $S\in\SV$ that has a positive 
 definite completion.
\item[\textbf{Output.}] The Cholesky factors $L$, $D$
 of $X=-\nabla f_*(S)$, \ie, the positive definite matrix 
 $X\in\SV$ that satisfies $\proj(X^{-1}) = S$. 
\item[\textbf{Algorithm.}]
Iterate over $j\in\{1,2,\ldots,n\}$ in reverse topological order.
For each $j$, compute $D_{jj}$ and the $j$th column of $L$ from
\BEQ \label{e-mf-compl-1}
 L_{\I{j}j} = -V_j^{-1} S_{\I{j}j}, \qquad 
 D_{jj} = (S_{jj} + S_{\I{j}j}^T L_{\I{j}j})^{-1},
\EEQ
and compute the update matrices
\BEQ \label{e-mf-compl-2}
  V_i = E_{\Ip{j}\I{i}}^T \left[\begin{array}{cc} 
  S_{jj} & S_{\I{j}j}^T \\ S_{\I{j}j} & V_j \end{array}\right] 
  E_{\Ip{j}\I{i}}, \quad i\in \ch(j).
\EEQ
\end{list}
\end{algdesc}

The update matrix $V_j$ can be discarded after cycle $j$, and 
$L_{\I{j}j}$ and $D_{jj}$ can overwrite $S_{\I{j}j}$ and $S_{jj}$ in a 
CCS data structure.

The cost of Algorithm~\ref{a-completion} is higher than that of 
Algorithm~\ref{a-projinv} because step~(\ref{e-mf-compl-1})
involves the solution of an equation with $V_j$ as coefficient matrix,
whereas~(\ref{e-mf-gradient-2}) only requires a multiplication.
It is therefore of interest to propagate a factorization of $V_j$
instead of the matrix itself, and to 
replace~(\ref{e-mf-compl-2}) by an efficient method for computing 
a factorization of $V_i$, given a factorization of $V_j$.
This idea can be implemented as follows.
We use a factorization of the form $V_j = R_jR_j^T$, with 
$R_j$ upper triangular of order $|\I{j}|$.  
We need to replace~(\ref{e-mf-compl-2}) with an efficient method for 
computing $R_i$ from $R_j$.
The matrix $E_{\Ip{j}\I{i}}$ can be partitioned as
\[
E_{\Ip{j}\I{i}}
= \left[\begin{array}{cc} 1 & 0 \\ 0 & E \end{array}\right],
\]
where $E = E_{\I{j}, \I{i}\setminus \{j\}}$.
(This follows from the fact that $j$ is the first element of 
$\Ip{j}$ by definition of $\Ip{j}$, and $j$ is also the first element 
of $\I{i}$ if $j$ is the parent of vertex $i$.  As a consequence, 
the 1,1 element of matrix $E_{\Ip{j}\I{i}}$ is equal to one).
From~(\ref{e-mf-compl-2}),
\[
 V_i 
 = E_{\Ip{j}\I{i}}^T \left[\begin{array}{cc}
   S_{jj} & S_{\I{j}j}^T \\ S_{\I{j}j} & V_j \end{array}\right] 
   E_{\Ip{j}\I{i}}
 = \left[\begin{array}{cc} a & b^T \\ b & C C^T\end{array}\right] 
\]
where 
\[
 a=S_{jj}, \qquad b= E^T S_{\I{j}j}, \qquad C = E^T R_j.
\]
The matrix $C$ is obtained from the upper triangular matrix $R_j$
by deleting the rows in $\Ip{j} \setminus \I{i}$.
It can be reduced to square upper triangular form by 
writing it as $C = R Q^T$ with $R$ upper triangular
and $Q$ orthogonal (a product of Householder transformations 
\cite{GoV:96}).  Then the triangular factor in $V_i = R_iR_i^T$  
is given by
\[
 R_i = \left[\begin{array}{cc}
   \left(a - \|R^{-1}b\|_2^2\right)^{1/2} 
  & \left(R^{-1}b\right)^T \\ 0 & R \end{array}\right].
\]
This is summarized below.

\begin{algdesc}{Matrix completion with factored update matrices}   
\label{a-completion-factored}
\begin{list}{}{}
\item[\textbf{Input.}] A matrix $S\in\SV$ that has a positive definite
 completion.
\item[\textbf{Output.}] The Cholesky factors $L$, $D$
 of $X=-\nabla f_*(S)$, \ie, of the positive definite matrix 
 $X\in\SV$ that satisfies $\proj(X^{-1}) = S$. 
\item[\textbf{Algorithm.}]
Iterate over $j\in\{1,2,\ldots,n\}$ in reverse topological order.
For each $j$, 
\BIT
\item compute $D_{jj}$ and the $j$th column of $L$ from
\BEQ \label{e-mf-compl-1-factored}
 L_{\I{j}j} = -R_j^{-T}R_j^{-1} S_{\I{j}j}, \qquad 
 D_{jj} = (S_{jj} + S_{\I{j}j}^T L_{\I{j}j})^{-1}
\EEQ
\item for $i\in \ch(j)$, compute a factorization
\BEQ \label{e-mfp-compl-3-factored}
 E_{\I{j}, \I{i}\setminus\{j\}}^T R_j = RQ^T
\EEQ
with $R$ upper triangular and $Q$ orthogonal, and compute
\BEQ \label{e-mf-compl-2-factored}
 R_i = \left[\begin{array}{cc}
   \left(S_{jj} - \|R^{-1}S_{\I{j}j} \|_2^2\right)^{1/2} 
  & \left(R^{-1}S_{\I{j}j}\right)^T \\ 0 & R \end{array}\right].
\EEQ
\mbox{}
\EIT
\end{list}
\end{algdesc}

The cost of step~(\ref{e-mf-compl-1-factored}) 
is order $|\I{j}|^2$  for the forward and back substitutions,
and the cost of~(\ref{e-mf-compl-2-factored}) is proportional
to $|\I{i}|^2$ (for the computation of $R^{-1}S_{\I{j}j}$).
The cost of the reduction to triangular form 
in~(\ref{e-mfp-compl-3-factored}) is difficult to quantify
because it depends on the number of rows in $R_j$ that are deleted in 
the multiplication $E^TR_j$ and on their positions.
However the total cost is usually much less than the cost of computing
$R_i$ from scratch, as shown by the experiments in the next section.

\subsection{Numerical results} 
\label{s-numerical-results} 
In this section, we give experimental results with
Algorithms~\ref{a-projinv}, \ref{a-completion}, and
\ref{a-completion-factored}.  The algorithms were implemented in
Python 2.7, using the Python library CVXOPT version 1.1.3
\cite{DaV:10} and its interfaces to LAPACK and BLAS\footnote{We link
  against the single-threaded reference implementations of BLAS and
  LAPACK in Ubuntu.} for the sparse and dense matrix computations.
Some critical code segments were implemented in C (such as the
Householder updates in Algorithm~\ref{a-completion-factored} and the
``extend-add'' operation $F := F + E_{JI} U E_{JI}^T$ and its
adjoint.)  The experiments were conducted on an Intel Q6600 CPU (2.4
GHz Core 2 Quad) computer with 4 GB memory, running Ubuntu 11.04.

\subsubsection{Band and arrow patterns}
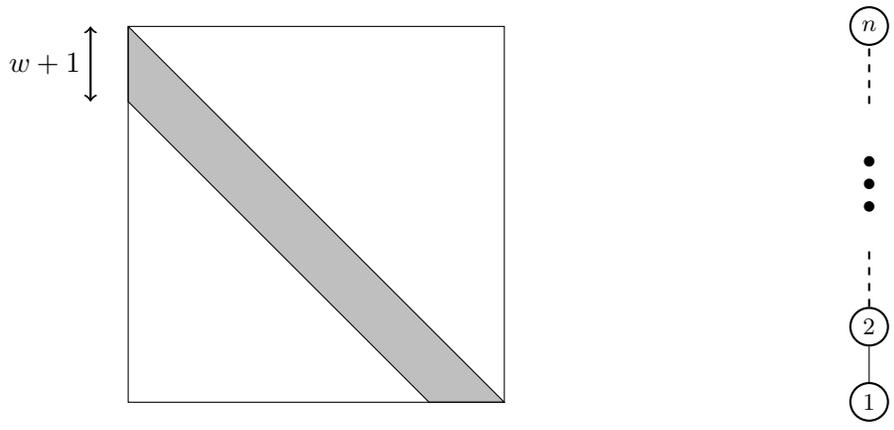
\begin{figure}
\hspace*{\fill}
\begin{minipage}{.5\linewidth}
\hspace*{\fill}
\begin{tikzpicture}[scale = 1.]
\draw [fill = gray!50] 
    (0,0) -- (0,-1) -- (4,-5) -- (5,-5) -- (0,0);
\draw (0,0) rectangle (5,-5);
\draw[<->, thick] (-.5,0) -- (-.5,-1);
\node[anchor = east] at (-.5,-.5) {$w+1$};
\end{tikzpicture}
\hspace*{\fill}
\end{minipage}
\hspace*{\fill}
\begin{minipage}{.45\linewidth}
\hspace*{\fill}
\begin{tikzpicture}[scale=1]
   \tikzset{VertexStyle/.style = {shape = circle, draw,  thick,
       minimum size = .5cm, inner sep = 2, fill=none, font=\footnotesize}}

      \node[VertexStyle](1) at (0,0){1};
   \node[VertexStyle](2) at (0,1){2};
   \node[VertexStyle, draw = none](3) at (0,2.3){};
   \node[VertexStyle, draw = none](4) at (0,3.7){};
   \node[VertexStyle](5) at (0,5){$n$};
   \foreach \x/\y in { 0/2.6, 0/2.9, 0/3.2 } \fill (\x,\y) circle (2pt); 

   \draw (1)--(2);
   \draw[dashed, thick] (2)--(3);
   \draw[dashed, thick] (4)--(5);
\end{tikzpicture}
\hspace*{\fill}
\end{minipage}
\hspace*{\fill}
\caption{Lower-triangular band pattern with bandwidth $w$ and the 
 corresponding elimination tree.  The vertices $1$, \ldots, $n-w$ have
 monotone degree $w$.  The vertices $k = n-w+1, \ldots,n$ have 
 monotone degree $n-k$.}
\label{f-band}
\end{figure}

\begin{figure}
\hspace*{\fill}
\begin{minipage}{.5\linewidth}
\hspace*{\fill}
\begin{tikzpicture}[scale = 1.]
\draw [fill = gray!50, draw = none] (0,-4) -- (0,-5) -- (5,-5) -- (4,-4) 
    -- (0,-4);
\draw (0,0) rectangle (5,-5);
\draw (0,0) -- (5,-5);
\draw (0,-4) -- (4,-4);
\draw[<->, thick] (-.5,-4) -- (-.5,-5);
\node[anchor = east] at (-.6,-4.5) {$w$};
\end{tikzpicture}
\hspace*{\fill}
\end{minipage}
\hspace*{\fill}
\begin{minipage}{.45\linewidth}
\hspace*{\fill}
\begin{tikzpicture}[scale=.8]
      \tikzset{VertexStyle/.style = {shape = circle, draw,  thick,
       minimum size = .5cm, inner sep = 2, fill=none, font=\small}}

      \node[VertexStyle](1) [label = 270: \footnotesize $1$] at (0,0){};
   \node[VertexStyle](2) [label = 270: \footnotesize $2$] at (2,0){};
   \node[VertexStyle](3) [label = 270: \footnotesize $3$] at (4,0){};
   \foreach \x/\y in {
       5.5/0, 6/0, 6.5/0
      }
      \fill (\x,\y) circle (2pt); 
   \node[VertexStyle](4) [label = 270: \footnotesize $n-w$] at (8,0){};
   \node[VertexStyle](5) [label = 0: \footnotesize $n-w+1 $] at (4,1.5){};
   \node(6) at (4,2.5){};
   \node(7) at (4,3.0){};
   \node[VertexStyle](8) [label = 0: \footnotesize $n-1$] at (4,4.0){};
   \node[VertexStyle](9) [label = 0: \footnotesize $n$] at (4,5.0){};

   \draw (1)--(5);
   \draw (2)--(5);
   \draw (3)--(5);
   \draw (4)--(5);
   \draw[dashed, thick] (5)--(6);
   \draw[dashed, thick] (7)--(8);
   \draw (8)--(9);
\end{tikzpicture}
\hspace*{\fill}
\end{minipage}
\hspace*{\fill}
\caption{Lower-triangular arrow pattern with width $w$ and 
 the corresponding elimination tree.  The vertices $1$, \ldots, $n-w$ have
 monotone degree $w$.  The vertices $k = n-w+1, \ldots,n$ have 
 monotone degree $n-k$.}
\label{f-arrow}
\end{figure}
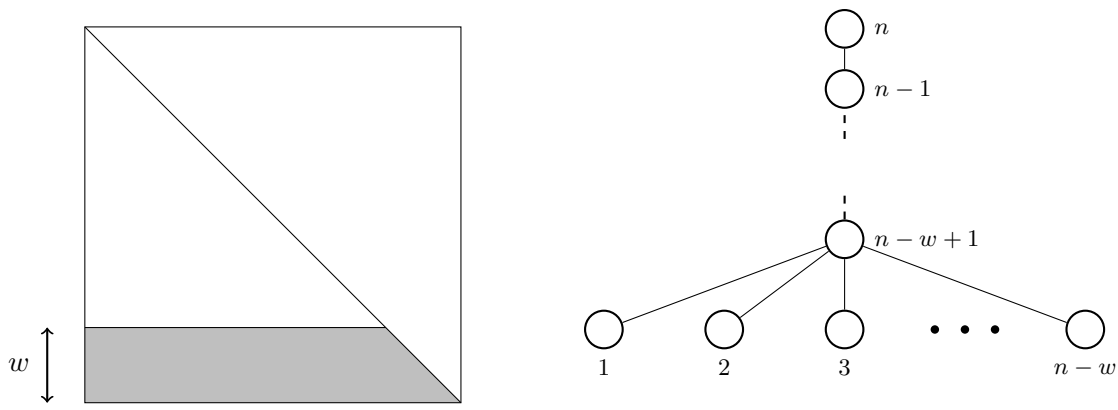

Band and arrow patterns are two basic sparsity patterns 
for which the complexity of the algorithms is easy to analyze.
We assume $n\gg w$, where $w$ is the bandwidth or blockwidth
(see Figures~\ref{f-band} and~\ref{f-arrow}).  
It is easy to see that the cost of a Cholesky factorization 
of a matrix with one of these two patterns is $O(nw^2)$.

Step~(\ref{e-mf-gradient-1}) of Algorithm~\ref{a-projinv}
involves matrix-vector multiplications of order $|I_j|$.
The complexity of the algorithm is dominated by the total cost of 
these products, \ie,  $O(nw^2)$.
This is similar to the cost of a Cholesky factorization.
Step~(\ref{e-mf-compl-1}) of Algorithm~\ref{a-completion} on the other
hand requires solving a dense positive definite system of order
$|I_j|$.  The total complexity is therefore $O(nw^3)$.

In Algorithm~\ref{a-completion-factored},
the cost of step~(\ref{e-mf-compl-1-factored})
is reduced to $O(w^2)$ per iteration.
For $j=n-w+1$, \ldots, $n$, we have $E=I$ in 
step~(\ref{e-mf-gradient-2}), so $C$ is upper triangular
and we only need to compute $R^{-1}b$.
For the other vertices in the elimination tree ($j=1,\ldots,n-w$),
$C= E^TR_j$ is the matrix $R_j$ with one row deleted:
the last row in the case of a band pattern, 
the first row in the case of an arrow pattern.
The cost of reducing $C$ to triangular form is therefore zero
in the case of an arrow pattern and $O(w^2)$ in the case of a band
pattern.  In either case, the total cost of 
Algorithm~\ref{a-completion-factored} is reduced to $O(nw^2)$.

\begin{figure}
\centering 
\psfrag{x00}[c][c]{${10^{1}}$}
\psfrag{x01}[c][c]{${10^{2}}$}
\psfrag{x02}[c][c]{${10^{3}}$}
\psfrag{y00}[r][r]{${10^{-2}}$}
\psfrag{y01}[r][r]{${10^{-1}}$}
\psfrag{y02}[r][r]{${10^{0}}$}
\psfrag{y03}[r][r]{${10^{1}}$}
\psfrag{y04}[r][r]{${10^{2}}$}
\psfrag{xlbl0}[c][b]{{\Large $w$}}
\psfrag{ylbl0}[c][t]{{\Large CPU time (seconds)}}
\psfrag{lgnd00www}[l][l]{Algorithm~\ref{a-projinv}}
\psfrag{lgnd01www}[l][l]{Algorithm~\ref{a-completion}}
\psfrag{lgnd02www}[l][l]{Algorithm~\ref{a-completion-factored}}
\resizebox{0.47\linewidth}{!}{\includegraphics[width=5in]{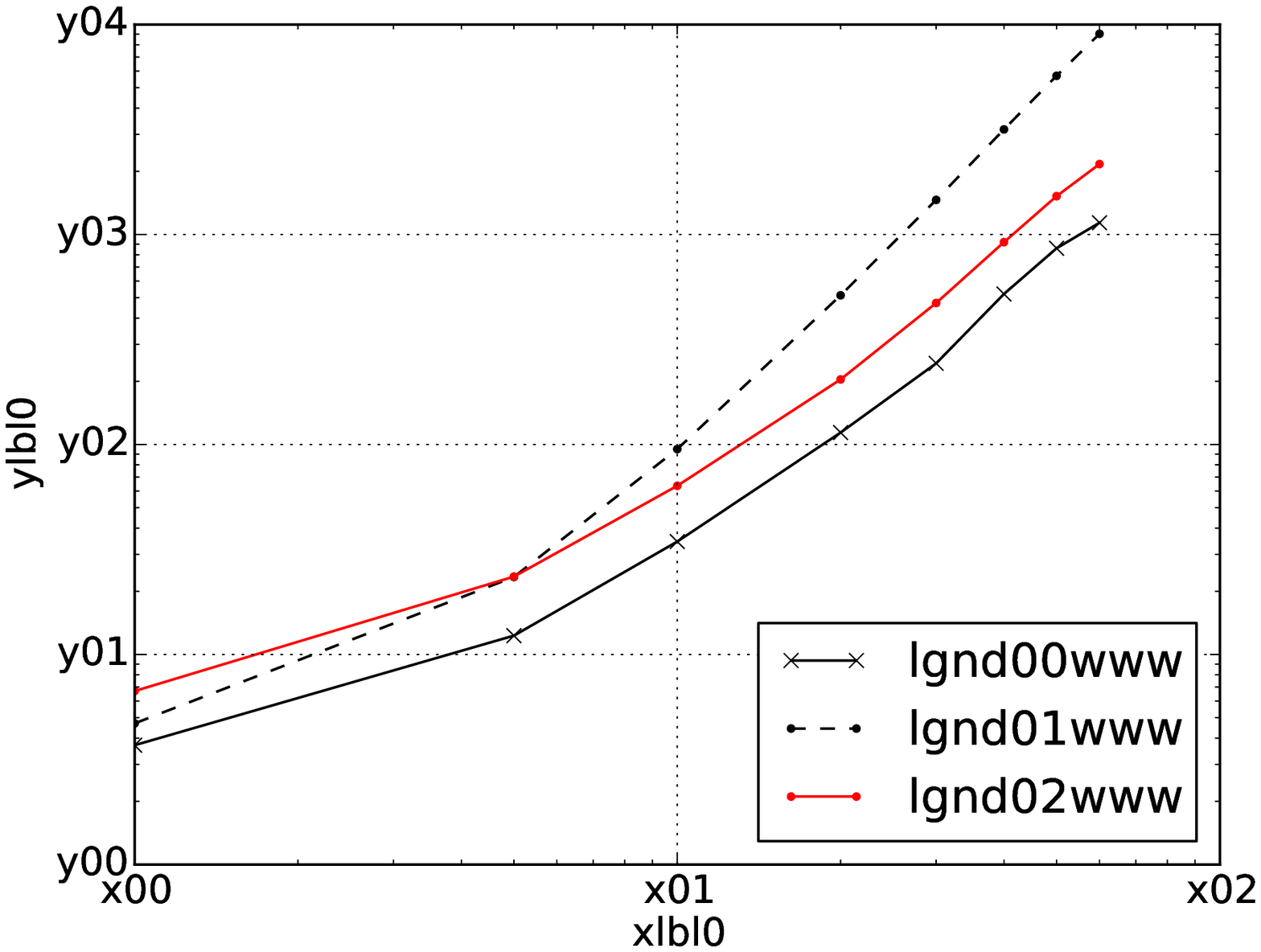}}\  
\psfrag{x00}[c][c]{${10^{1}}$}
\psfrag{x01}[c][c]{${10^{2}}$}
\psfrag{x02}[c][c]{${10^{3}}$}
\psfrag{y00}[r][r]{${10^{-2}}$}
\psfrag{y01}[r][r]{${10^{-1}}$}
\psfrag{y02}[r][r]{${10^{0}}$}
\psfrag{y03}[r][r]{${10^{1}}$}
\psfrag{y04}[r][r]{${10^{2}}$}
\psfrag{xlbl0}[c][b]{{\Large $w$}}
\psfrag{ylbl0}[c][t]{{\Large CPU time (seconds)}}
\psfrag{lgnd00www}[l][l]{Algorithm~\ref{a-projinv}}
\psfrag{lgnd01www}[l][l]{Algorithm~\ref{a-completion}}
\psfrag{lgnd02www}[l][l]{Algorithm~\ref{a-completion-factored}}
\resizebox{0.47\linewidth}{!}{\includegraphics[width=5in]{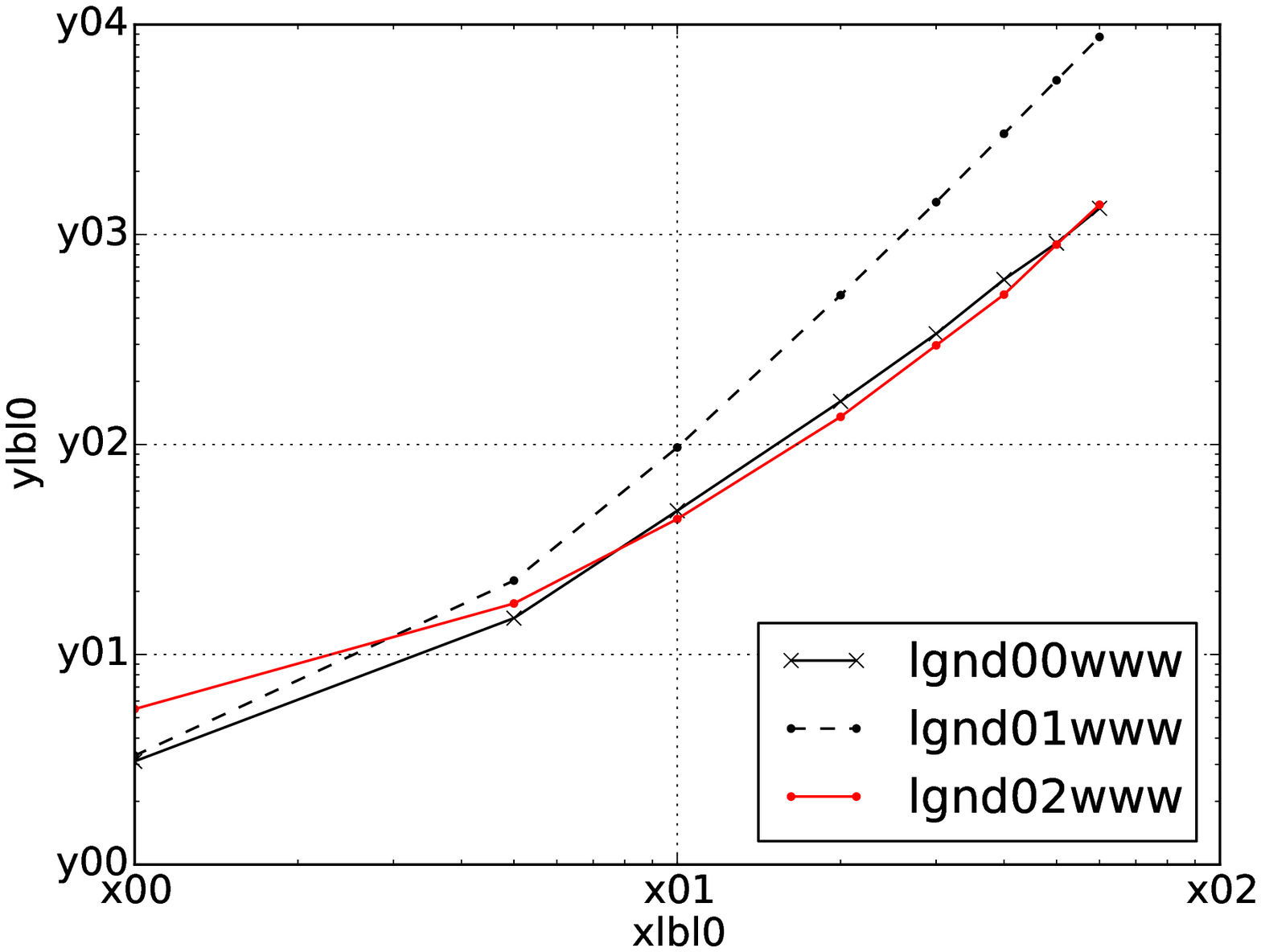}}\caption{Cost of computing the primal gradient 
(Algorithm~\ref{a-projinv}) and the dual gradient 
(Algorithms~\ref{a-completion} and~\ref{a-completion-factored}) 
for (a) band patterns and (b) arrow patterns of 
order $n=2000$, as a function of the pattern width $w$.}  
\label{f-band-arrow-ex}
\end{figure}

In Figure~\ref{f-band-arrow-ex}, we compare the CPU times
of the three algorithms as a function of the width $w$.
As can be seen, the cost of computing the dual gradient using
Algorithm~\ref{a-completion-factored} is comparable to the cost of the
primal gradient using Algorithm~\ref{a-projinv}, and the cost of 
the two algorithms grows roughly as $w^2$ for fixed $n$. 
Notice the small gap between the cost of
Algorithms~\ref{a-completion-factored} 
and~\ref{a-projinv} for band patterns. This gap reflects the
cost of reducing $C$ to triangular form
in Algorithm~\ref{a-completion-factored}.

\subsubsection{General sparse patterns}
\begin{figure}
  \centering
\psfrag{x00}[c][c]{${10^{2}}$}
\psfrag{x01}[c][c]{${10^{3}}$}
\psfrag{x02}[c][c]{${10^{4}}$}
\psfrag{x03}[c][c]{${10^{5}}$}
\psfrag{x04}[c][c]{${10^{6}}$}
\psfrag{y00}[r][r]{${10^{-5}}$}
\psfrag{y01}[r][r]{${10^{-4}}$}
\psfrag{y02}[r][r]{${10^{-3}}$}
\psfrag{y03}[r][r]{${10^{-2}}$}
\psfrag{y04}[r][r]{${10^{-1}}$}
\psfrag{y05}[r][r]{${10^{0}}$}
\psfrag{xlbl0}[c][b]{{\Large $n$}}
\psfrag{ylbl0}[c][t]{{\Large Density}}
\resizebox{0.47\linewidth}{!}{\includegraphics[width=5in]{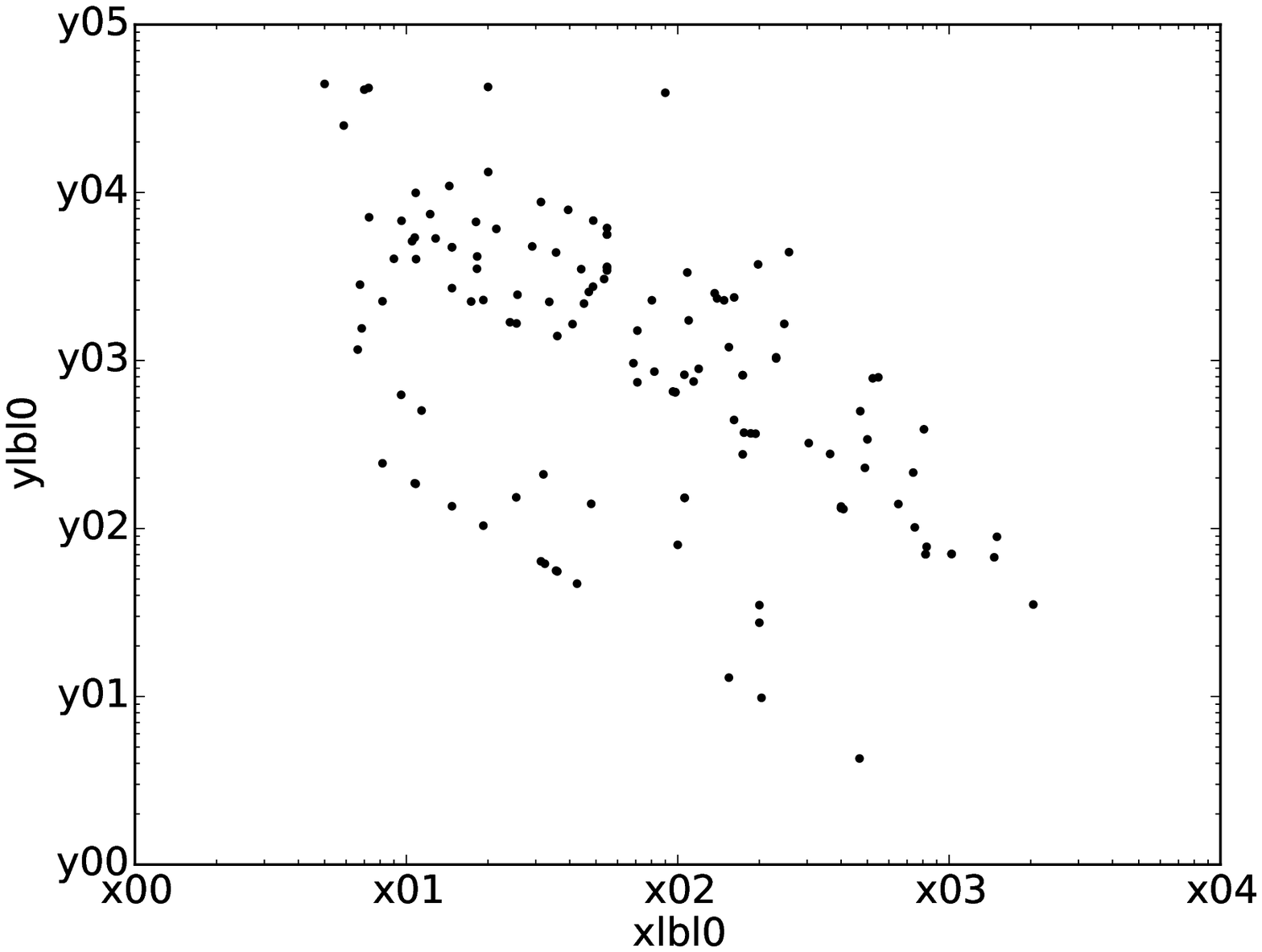}}  \
\psfrag{x00}[c][c]{${10^{2}}$}
\psfrag{x01}[c][c]{${10^{3}}$}
\psfrag{x02}[c][c]{${10^{4}}$}
\psfrag{x03}[c][c]{${10^{5}}$}
\psfrag{x04}[c][c]{${10^{6}}$}
\psfrag{y00}[r][r]{${10^{2}}$}
\psfrag{y01}[r][r]{${10^{3}}$}
\psfrag{y02}[r][r]{${10^{4}}$}
\psfrag{y03}[r][r]{${10^{5}}$}
\psfrag{y04}[r][r]{${10^{6}}$}
\psfrag{y05}[r][r]{${10^{7}}$}
\psfrag{y06}[r][r]{${10^{8}}$}
\psfrag{xlbl0}[c][b]{{\Large $n$}}
\psfrag{ylbl0}[c][t]{{\Large Number of noneros in $L$}}
\resizebox{0.47\linewidth}{!}{\includegraphics[width=5in]{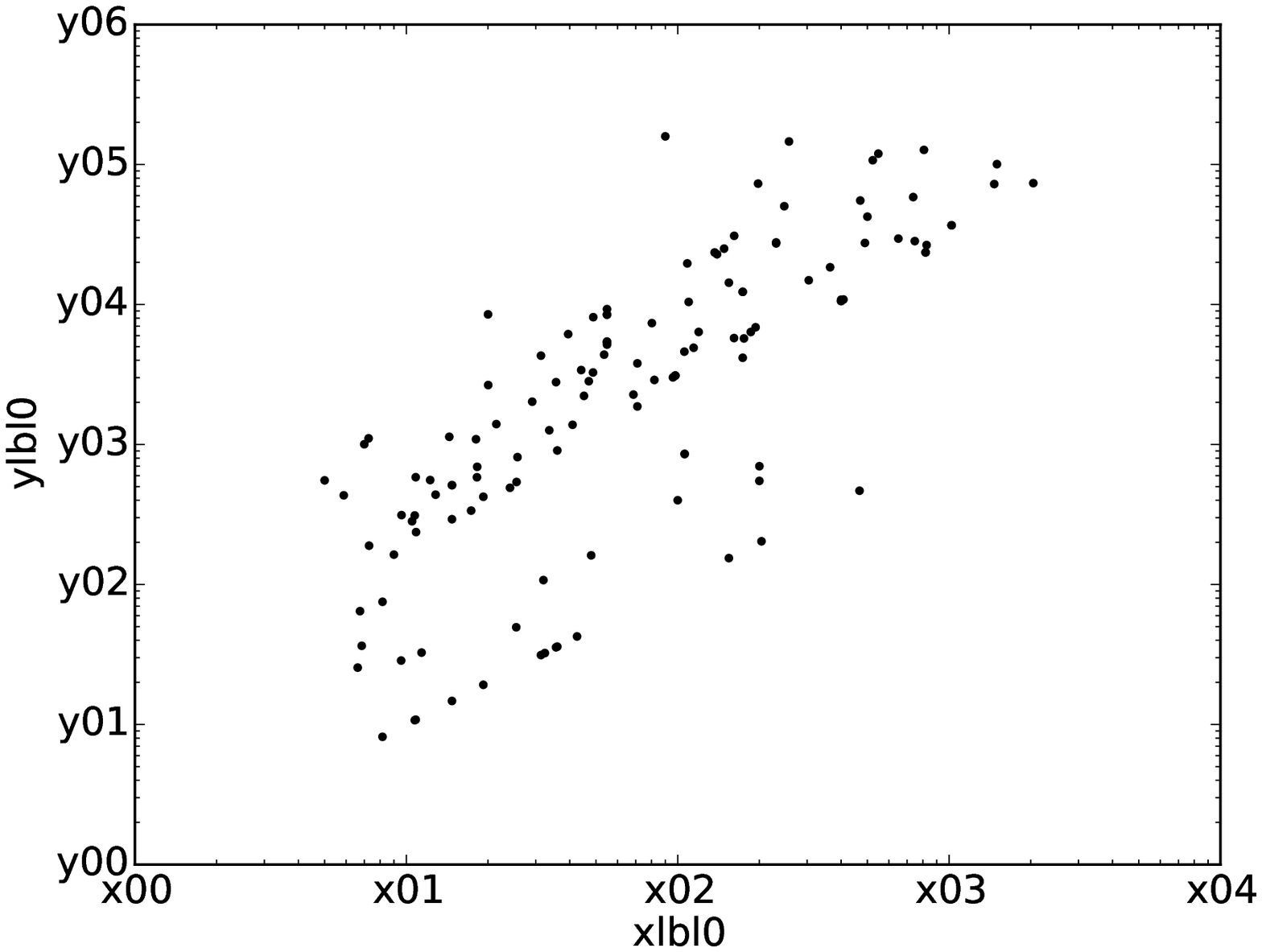}}  \caption{Scatter plot of $n$ versus (a) the density
    ($|V|/(n(n+1)/2)$) and (b) the number of nonzeros in $L$ for the
    128 test problems. Each dot represents a problem from the
    University of Florida Sparse Matrix collection.}
  \label{fig-ex_gradient_ufsmc_stat}
\end{figure}

In the second experiment, we use a benchmark set of large symmetric
sparsity patterns from the University of Florida Sparse Matrix
Collection \cite{Dav:09b}.  The AMD ordering was used to compute
filled patterns.  To prevent out-of-core computations, we restrict the
experiment to matrices for which the filled pattern occupies less than 250 MB of
memory. The set of test problems includes 128 sparsity patterns, with
$n$ ranging from 500 to 204316 and with $|V|$ between 817 and
15,894,180. A scatter plot of the number of nonzeros and the density of
the test problems versus the dimension $n$ is shown in 
Figure~\ref{fig-ex_gradient_ufsmc_stat}.

\begin{figure}
  \centering
\psfrag{x00}[c][c]{${10^{-3}}$}
\psfrag{x01}[c][c]{${10^{-2}}$}
\psfrag{x02}[c][c]{${10^{-1}}$}
\psfrag{x03}[c][c]{${10^{0}}$}
\psfrag{x04}[c][c]{${10^{1}}$}
\psfrag{x05}[c][c]{${10^{2}}$}
\psfrag{x06}[c][c]{${10^{3}}$}
\psfrag{y00}[r][r]{${10^{-3}}$}
\psfrag{y01}[r][r]{${10^{-2}}$}
\psfrag{y02}[r][r]{${10^{-1}}$}
\psfrag{y03}[r][r]{${10^{0}}$}
\psfrag{y04}[r][r]{${10^{1}}$}
\psfrag{y05}[r][r]{${10^{2}}$}
\psfrag{y06}[r][r]{${10^{3}}$}
\psfrag{xlbl0}[c][b]{{\Large Algorithm~\ref{a-projinv}: projected
 inverse (sec.)}}
\psfrag{ylbl0}[c][t]{{\Large Algorithm~\ref{a-completion-factored}:
 fast completion (sec.)}}
\resizebox{0.47\linewidth}{!}{\includegraphics[width=5in]{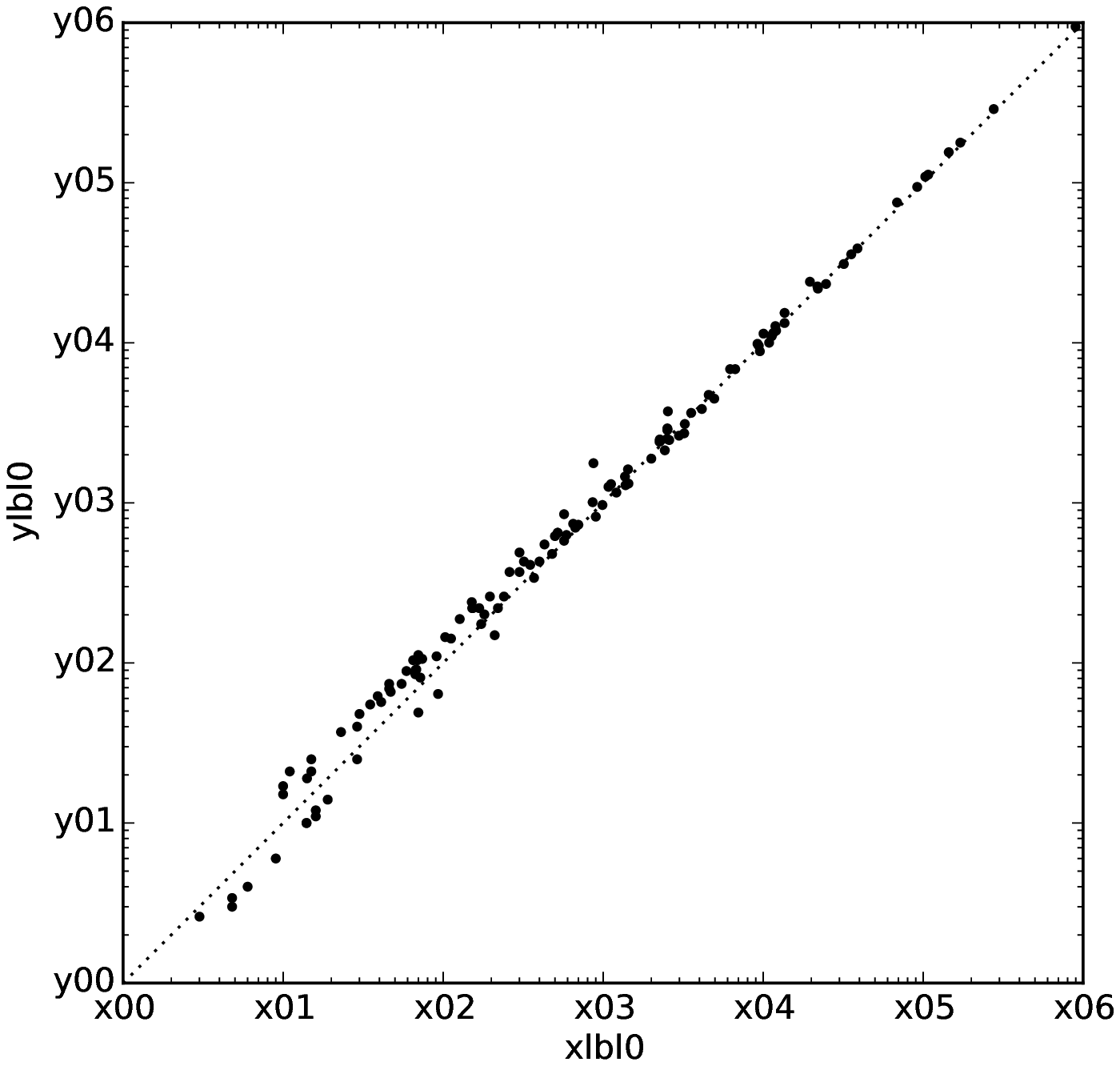}}
  \ 
\psfrag{x00}[c][c]{${10^{-3}}$}
\psfrag{x01}[c][c]{${10^{-2}}$}
\psfrag{x02}[c][c]{${10^{-1}}$}
\psfrag{x03}[c][c]{${10^{0}}$}
\psfrag{x04}[c][c]{${10^{1}}$}
\psfrag{x05}[c][c]{${10^{2}}$}
\psfrag{x06}[c][c]{${10^{3}}$}
\psfrag{x07}[c][c]{${10^{4}}$}
\psfrag{y00}[r][r]{${10^{-3}}$}
\psfrag{y01}[r][r]{${10^{-2}}$}
\psfrag{y02}[r][r]{${10^{-1}}$}
\psfrag{y03}[r][r]{${10^{0}}$}
\psfrag{y04}[r][r]{${10^{1}}$}
\psfrag{y05}[r][r]{${10^{2}}$}
\psfrag{y06}[r][r]{${10^{3}}$}
\psfrag{y07}[r][r]{${10^{4}}$}
\psfrag{xlbl0}[c][b]{{\Large Algorithm~\ref{a-completion}: 
 completion (sec.)}}
\psfrag{ylbl0}[c][t]{{\Large Algorithm~\ref{a-completion-factored}:
fast completion (seconds)}}
\resizebox{0.47\linewidth}{!}{\includegraphics[width=5in]{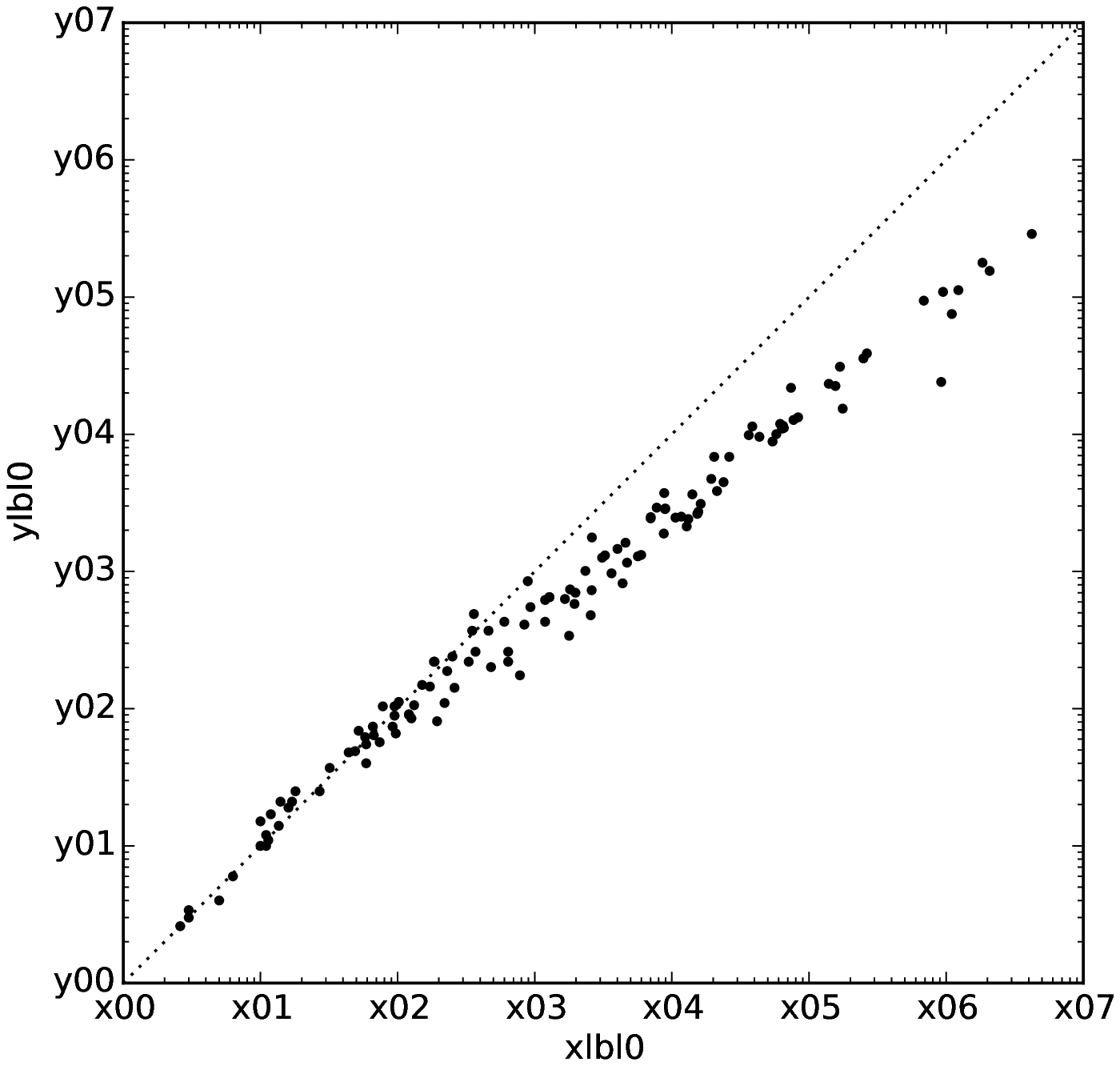}}
  \caption{CPU times for primal and dual gradient algorithms  for 
  a test set of 128 sparsity patterns: (a) shows the cost
  of Algorithm~\ref{a-projinv} (primal gradient or projected inverse) 
  versus the cost of Algorithm~\ref{a-completion-factored} 
  (dual gradient or matrix completion) for different
  sparsity patterns; (b) compares the cost of computing the dual
  gradient via Algorithms~\ref{a-completion-factored} 
  and~\ref{a-completion} 
  \ie, with and without the update technique described in
  Section~\ref{s-dual-gradient}.}
  \label{fig-primal-dual-gradients}
\end{figure}

Figure~\ref{fig-primal-dual-gradients} shows the CPU times 
for Algorithm~\ref{a-projinv} (primal gradient or projected inverse)
and~\ref{a-completion-factored} (dual gradient or completion), 
and for Algorithms~\ref{a-completion-factored} and~\ref{a-completion} 
(completion with and without Householder updates, respectively).
Each dot represents one of the sparsity patterns in the test set.  
The results indicate that in practice, on this set of realistic sparsity 
patterns, the costs of computing the primal and dual gradients 
are comparable. 

\section{Hessian} \label{s-hess}
The Hessian $\mathcal H = \nabla^2f(X)$ of $f$ at $X$ is defined as 
\[
\mathcal H(Y) =
 \proj(X^{-1}YX^{-1}) = -\left.\frac{d}{dt} \proj(X+tY)^{-1} \right|_{t=0}.
\]
A method for evaluating $\mathcal H(Y)$ can therefore be found by 
differentiating the formulas for evaluating the gradient 
$-\proj(X^{-1})$.  
As we have seen, $\proj(X^{-1})$ is obtained in two stages. 
First the Cholesky factors $L$, $D$ of $X$ are computed, column by column,
following a topological ordering of the elimination tree 
(Algorithm~\ref{a-chol}).  
Then the projected inverse $\proj(X^{-1})$ is computed from $L$ and $D$, 
column by column, in reverse topological order 
(Algorithm~\ref{a-projinv}).  
Linearizing the two algorithms will provide an algorithm for 
$\mathcal H(Y)$.  We give the details in Section~\ref{s-hessian}.

We also consider the problem of evaluating $Y = \mathcal H^{-1}(T)$, 
\ie, solving the linear equation 
\[
 \proj(X^{-1}YX^{-1}) = T
\]
for $Y\in\SV$, given $T\in\SV$.   An algorithm for this problem can
be formulated by inverting the calculation of $\mathcal H(Y)$ or,
alternatively, by linearizing the algorithms for matrix 
completion (Algorithms~\ref{a-completion} and~\ref{a-completion-factored}) 
and the Cholesky product (Algorithm~\ref{a-chol-prod}); 
see~Section~\ref{s-inv-hessian}.

In applications, it is often useful to know a factorization of the 
Hessian as a composition of a mapping and its adjoint,
\[
 \mathcal H(Y) = \mathcal R^\mathrm{adj}(\mathcal R(Y)).
\]
Such a factorization is discussed in Section~\ref{s-hessian-fact}.

\subsection{Linearized Cholesky factorization and matrix completion}
Let $L(t)$, $D(t)$ be the matrices in the factorization 
$X(t) = X+tY = L(t) D(t) L(t)^T$  and let $U_i(t)$ be the $i$th update 
matrix in the multifrontal factorization algorithm for $X(t)$, \ie,
\[
   U_i(t) = -\sum_{k\in T_i} D_{kk}(t) L_{\I{k}k}(t) L_{\I{k}k}(t)^T.
\]
We denote by $L'$, $D'$, $U_i'$ the derivatives of $L(t)$, $D(t)$, 
$U_i(t)$ at $t=0$.   These derivatives can be found by linearizing the 
equation~(\ref{e-mf-chol-0}) with $X$ replaced by $X+tY$,
\BEAS
\lefteqn{
\left[\begin{array}{cc} 
 X_{jj} & X_{\I{j}j}^T \\ X_{\I{j}j} & 0 \end{array}\right]
+ t \left[\begin{array}{cc} 
 Y_{jj} & Y_{\I{j}j}^T \\ Y_{\I{j}j} & 0 \end{array}\right]
 + \sum_{i\in\ch(j)} E_{\Ip{j}\I{i}} U_i(t) E_{\Ip{j}\I{i}}^T } \\
 & = & 
 D_{jj}(t) \left[\begin{array}{c} 1 \\ L_{\I{j}j}(t) \end{array}\right]
 \left[\begin{array}{c} 1 \\ L_{\I{j}j}(t) \end{array}\right]^T
 + \left[\begin{array}{cc} 0 & 0 \\ 0 & U_j(t) \end{array}\right].
\EEAS
Taking the derivatives of the left- and right-hand sides at $t=0$ gives
\BEQ \label{e-mf-chol-lin}
\left[\begin{array}{cc} 
 Y_{jj} & Y_{\I{j}j}^T \\ Y_{\I{j}j} & 0 \end{array}\right]
 + \sum_{i\in\ch(j)} E_{\Ip{j}\I{i}} U_i' E_{\Ip{j}\I{i}}^T 
  = 
  \left[\begin{array}{cc} 1 & 0 \\ L_{\I{j}j} & I \end{array}\right]
  \left[\begin{array}{cc} D_{jj}' & (D_{jj}L_{\I{j}j}')^T \\ 
   D_{jj}L_{\I{j}j}'  & U_j'\end{array}\right]
  \left[\begin{array}{cc} 1 & L_{\I{j}j}^T \\ 0 & I \end{array}\right].
\EEQ
This will be the key equation for computing the linearized Cholesky 
factors $D'$, $L'$ from $Y$, and conversely, the linearized Cholesky 
product $Y$ from the linearized factors $D'$, $L'$.

Similarly, we define $Z(t) = (X+tY)^{-1}$, $S(t) = \proj(Z(t))$, and
\[
 V_i(t) = S_{\I{i}\I{i}}(t).
\]
We write the derivatives of $S(t)$ and $V_i(t)$ at $t=0$ as $-T$ and 
$-V_i'$. 
Substituting $S(t)$, $L(t)$, $D(t)$ for $S$, $L$, $D$ 
in~(\ref{e-ldl-inv}) gives
\[
 \left[\begin{array}{cc} S_{jj}(t) & S_{\I{j}j}(t)^T \\ 
 S_{\I{j}j}(t) & S_{\I{j}\I{j}}(t)
 \end{array}\right]
 \left[\begin{array}{c} 1 \\ L_{\I{j}j}(t) \end{array} \right] = 
 \left[\begin{array}{c} 1/D_{jj}(t) \\ 0 \end{array}\right] 
\]
and differentiating with respect to $t$ gives
\[
-\left[\begin{array}{cc}
 T_{jj}  & T_{\I{j}j}^T \\ T_{\I{j}j} & V_j' \end{array}\right]
\left[\begin{array}{c} 1 \\ L_{\I{j}j} \end{array}\right] + 
 \left[\begin{array}{cc} S_{jj} & S_{\I{j}j}^T \\ 
 S_{\I{j}j} & S_{\I{j}\I{j}} \end{array}\right] 
 \left[\begin{array}{c} 0 \\ L_{\I{j}j}'\end{array} \right]
 = \left[\begin{array}{c} -D_{jj}'/D_{jj}^2  \\ 0 \end{array}\right].
\]
Using $S_{\I{j}j} = -S_{\I{j}\I{j}} L_{\I{j}j}$ 
(from~(\ref{e-grad})) this can be written in a 
more symmetric form as 
\[
\left[\begin{array}{cc}
 T_{jj}  & T_{\I{j}j}^T \\ T_{\I{j}j} & V_j' \end{array}\right]
\left[\begin{array}{cc} 1 & 0 \\ L_{\I{j}j} & I \end{array}\right]
 = \left[\begin{array}{cc} 1 & -L_{\I{j}j}^T  \\ 0 & I \end{array}\right]
  \left[\begin{array}{cc} D_{jj}'/D_{jj}^2  & (S_{\I{j}\I{j}}L_{\I{j}j}')^T  \\ 
   S_{\I{j}\I{j}}L_{\I{j}j}' & V_j' \end{array}\right], 
\]
\ie,
\BEQ \label{e-compl-lin}
\left[\begin{array}{cc} 1 & L_{\I{j}j}^T  \\ 0 & I \end{array}\right]
\left[\begin{array}{cc}
 T_{jj}  & T_{\I{j}j}^T \\ T_{\I{j}j}^T & V_j' \end{array}\right]
\left[\begin{array}{cc} 1 & 0 \\ L_{\I{j}j} & I \end{array}\right]
 = \left[\begin{array}{cc} D_{jj}'/D_{jj}^2  & 
  (S_{\I{j}\I{j}}L_{\I{j}j}')^T  \\ 
   S_{\I{j}\I{j}}L_{\I{j}j}' & V_j' \end{array}\right]. 
\EEQ
This equation allows us to compute $T$ given the linearized factors
$D'$, $L'$, and conversely, compute $D'$, $L'$ given $T$.

\subsection{Hessian} \label{s-hessian}
The algorithm for computing $T = \mathcal H(Y)$ first computes
$L'$, $D'$ by the linearized Cholesky factorization, \ie, 
from~(\ref{e-mf-chol-lin}), and then $T$ from $L'$, $D'$ 
by the linearized projected inverse algorithm, \ie, from~(\ref{e-compl-lin}).
To simplify the notation, we define two matrices $K$, $M\in\SV$ 
as
\BEQ\label{e-mf-hess-vars}
 K_{jj} = D_{jj}', \qquad K_{\I{j}j} = L_{\I{j}j}'D_{jj}, \qquad 
 M_{jj} = \frac{D_{jj}'}{D_{jj}^{2}},  \qquad
 M_{\I{j}j} = S_{\I{j}\I{j}}L_{\I{j}j}', 
\EEQ
for $j=1,\ldots,n$.

\newpage
\begin{algdesc}{Hessian evaluation}  \label{a-hessian}
\begin{list}{}{}
\item[\textbf{Input.}] A matrix $Y\in\symm^n_V$, the Cholesky factors 
$L$, $D$ of a positive definite matrix $X\in\SV$,
and the projected inverse $S=\proj(X^{-1})$.
\item[\textbf{Output.}] The matrix 
$T= \mathcal H(Y) = \proj(X^{-1}YX^{-1})$ where $\mathcal H$
is the Hessian of $f$ at $X$. 
\item[\textbf{Algorithm.}] \mbox{}
\begin{enumerate}
\item Iterate over $j\in\{1,2,\ldots,n\}$ in topological order.
 For each $j$, calculate $D_{jj}'$, the $j$th column of $K$, and 
 the update matrix $U_j'$ via
\BEAS
 \lefteqn{
 \left[\begin{array}{cc} K_{jj} & K_{\I{j}j}^T \\ K_{\I{j}j} & U_j' 
 \end{array}\right] } \\
 & = & \left[\begin{array}{cc} 1 & 0 \\ -L_{\I{j}j} & I \end{array}\right] 
 \left(\left[\begin{array}{cc}
     Y_{jj} & Y_{\I{j}j}^T \\ Y_{\I{j}j} & 0 \end{array}\right]
   + \sum_{i\in\ch(j)} E_{\Ip{j}\I{i}} U_i' E_{\Ip{j}\I{i}}^T \right)
  \left[\begin{array}{cc} 1 & -L_{\I{j}j}^T \\ 0 & I \end{array}\right].
\EEAS

\item For $j \in \{1,2,\ldots,n\}$, compute column $j$ of $M$ via
\[
 M_{jj} = \frac{K_{jj}}{D_{jj}^2}, \qquad 
 M_{\I{j}j} = \frac{1}{D_{jj}} S_{\I{j}\I{j}} K_{\I{j}j}. 
\]

\item Iterate over $j\in \{1,2,\ldots,n\}$ in reverse topological order.
 For each $j$, calculate $T_{jj}$ and $T_{\I{j}j}$ from 
\[
 \left[\begin{array}{cc}
  T_{jj} & T_{\I{j}j}^T \\ T_{\I{j}j} & V_j' \end{array}\right] =
  \left[\begin{array}{cc} 1 & -L_{\I{j}j}^T  \\ 0 & I \end{array}\right]
  \left[\begin{array}{cc}
  M_{jj}  &  M_{\I{j}j}^T \\ M_{\I{j}j} & V_j' \end{array}\right] 
  \left[\begin{array}{cc} 1 & 0 \\ -L_{\I{j}j} & I \end{array}\right] 
 \]
and the update matrices $V_i'$ for the children of vertex $j$ via
 \[
  V_i' = E_{\Ip{j}\I{i}}^T 
 \left[\begin{array}{cc} T_{jj} & T_{\I{j}j}^T \\ 
  T_{\I{j}j} & V_j' \end{array}\right] E_{\Ip{j}\I{i}}, \quad i\in \ch(j). 
 \]
\end{enumerate}
\end{list}
\end{algdesc}

The vertices $j$ in step~2 can be ordered in any order.
However, by defining $V_j = S_{\I{j}\I{j}}$
as in Algorithm~\ref{a-projinv} and using a reverse topological ordering,
we can avoid having to extract $S_{\I{j}\I{j}}$ from the CCS 
structure of $S$. 
In the modified algorithm, the second step is replaced by 
\[
 M_{jj} = \frac{K_{jj}}{D_{jj}^2}, \qquad
 M_{\I{j}j} = \frac{1}{D_{jj}} V_j K_{\I{j}j},  \qquad
 V_i = E_{\Ip{j}\I{i}}^T \left[\begin{array}{cc} 
  S_{jj} & S_{\I{j}j}^T \\ S_{\I{j}j} & V_j \end{array}\right] 
   E_{\Ip{j}\I{i}}, \quad i \in \ch(j),
\]
for $j \in \{1,2,\ldots,n\}$ in a reverse topological order.

\subsection{Inverse Hessian} \label{s-inv-hessian}
To evaluate $Y = \mathcal H^{-1}(T)$, we use the 
equation~(\ref{e-mf-chol-lin}) to compute the linearized Cholesky
factors $D'$, $L'$ from $T$, and the equation~(\ref{e-compl-lin}) to 
compute $Y$ from $D'$, $L'$.
We use the same notation~(\ref{e-mf-hess-vars}) 
as in the previous section.

\begin{algdesc}{Inverse Hessian evaluation} \label{a-inv-hessian}
\begin{list}{}{}
\item[\textbf{Input.}]  A matrix $T\in\SV$, the
Cholesky $L$, $D$ of a positive definite matrix $X \in\SV$, and 
the projected inverse $S = \proj(X^{-1})$. 
\item[\textbf{Output.}] The matrix $Y=\mathcal H^{-1}(T)$, \ie,
 the solution $Y\in\SV$ of the equation $\mathcal H(Y) = 
 \proj(X^{-1}YX^{-1}) = T$.
\item[\textbf{Algorithm.}] \mbox{}
\begin{enumerate}
\item Iterate over $j\in\{1,2,\ldots,n\}$ in reverse topological order.
 For each $j$, calculate the $j$th column of $M$ from
\[
\left[\begin{array}{cc} M_{jj} & M_{\I{j}j}^T \\ M_{\I{j}j} & V_j' 
 \end{array}\right]
 = \left[\begin{array}{cc} 1 & L_{\I{j}j}^T \\ 0 & I \end{array}\right]
 \left[\begin{array}{cc} T_{jj} & T_{\I{j}j}^T \\ T_{\I{j}j} & V_j'
 \end{array}\right] 
 \left[\begin{array}{cc} 1 & 0 \\ L_{\I{j}j} & I \end{array}\right]
\]
and the update matrices $V_i'$ for the children of vertex $j$ via
\[
  V_i' = E_{\Ip{j}\I{i}}^T \left[\begin{array}{cc} 
 T_{jj} & T_{\I{j}j}^T \\ T_{\I{j}j} & V_j' \end{array}\right] 
  E_{\Ip{j}\I{i}}, \quad i\in \ch(j).
\]

\item For $j \in \{1,2,\ldots,n\}$, compute column $j$ of $K$ via
\[
 K_{jj} = D_{jj}^2M_{jj}, \qquad K_{\I{j}j} = D_{jj} S_{\I{j}\I{j}}^{-1}  
  M_{\I{j}j}.
\]

\item Iterate over $j\in\{1,2,\ldots,n\}$ in topological order.
 For each $j$, compute $U_j'$ and the $j$th column of $Y$ from
\[
\left[\begin{array}{cc}
 Y_{jj} & Y_{\I{j}j}^T \\ Y_{\I{j}j} & -U_j' \end{array}\right]
= \left[\begin{array}{cc} 1 & 0 \\ L_{\I{j}j} & I \end{array}\right]
 \left[\begin{array}{cc}
  K_{jj} & K_{\I{j}}^T \\ K_{\I{j}j} & 0 \end{array}\right]
  \left[\begin{array}{cc} 1 & L_{\I{j}j}^T \\ 0 & I \end{array}\right]
 - \sum_{i\in\ch(j)} E_{\Ip{j}\I{i}} U_i' E_{\Ip{j}\I{i}}^T.
\]
\end{enumerate}
\end{list}
\end{algdesc}
A improvement of step~2 is to use a factorization 
$S_{\I{j}\I{j}} = -R_j R_j^T$  
and update the matrices $R_j$ recursively, following a reverse 
topological order, as discussed at the end of~Section~\ref{s-grad}.

\subsection{Hessian factor} \label{s-hessian-fact}
Step~1 in the Hessian evaluation algorithm (Algorithm~\ref{a-hessian})
is a linear mapping that transforms $Y$ to $K$. 
Step~3 is a linear mapping that transforms $M$ to $T$.  
It is interesting to note that these two mappings are adjoints.  
Step~2 implements a self-adjoint and positive definite mapping, which 
transforms $K$ to $M$.
Factoring the positive definite mapping in step~2 provides
a factorization 
\[ 
  \mathcal H(Y) = \mathcal R^\mathrm{adj}(\mathcal R(Y)),
\]
with $\mathcal R$ a linear mapping from $\symm^n_V$ to $\symm^n_V$.
The factorization of the mapping in step~2 can be implemented by
defining a factorization
$S_{\I{j}\I{j}} = R_jR_j^T$
with $R_j$ upper triangular for each vertex $j$.
Although it is impractical to pre-compute and store the matrices $R_j$ 
for each vertex, they can be efficiently computed recursively in a reverse
topological order, as in Algorithm~\ref{a-completion-factored}.
For the sake of clarity, we omit the details in the following algorithms.

The algorithm for evaluating $\mathcal R$ consists 
of step~1 in Algorithm~\ref{a-hessian} and one half of step~2.

\begin{algdesc}{Evaluation of Hessian factor} \label{a-hessian-factor}
\begin{list}{}{}
\item[\textbf{Input.}]
A matrix $Y\in\symm^n_V$,
the Cholesky factors $L$, $D$ of a positive definite matrix 
$X\in\SV$, and the projected inverse $S = \proj(X^{-1})$. 
\item[\textbf{Output.}] The matrix $W = \mathcal R(Y)$
where $\mathcal H = \mathcal R^\mathrm{adj}\circ \mathcal R$ is the 
Hessian of $f$ at $X$.
\item[\textbf{Algorithm.}] \mbox{}
\begin{enumerate}
\item Iterate over $j\in\{1,2,\ldots,n\}$ in topological order.
 For each $j$, calculate the $j$th column of $K$, and 
 the update matrix $U_j'$ via 
\BEAS
\lefteqn{
 \left[\begin{array}{cc} K_{jj} & K_{\I{j}j}^T \\ K_{\I{j}j} & U_j' 
 \end{array}\right] } \\
 & = & \left[\begin{array}{cc} 1 & 0 \\ -L_{\I{j}j} & I \end{array}\right] 
  \left( \left[\begin{array}{cc}
     Y_{jj} & Y_{\I{j}j}^T \\ Y_{\I{j}j} & 0 \end{array}\right]
   + \sum_{i\in\ch(j)} E_{\Ip{j}\I{i}} U_i' E_{\Ip{j}\I{i}}^T\right)
  \left[\begin{array}{cc} 1 & -L_{\I{j}j}^T \\ 0 & I \end{array}\right].
\EEAS

\item For all $j \in \{1,2,\ldots,n\}$, compute column $j$ of $W$ via
\[
 W_{jj} = \frac{K_{jj}}{D_{jj}}, \qquad 
 W_{\I{j}j} = \frac{1}{\sqrt D_{jj}} R_j^T  K_{\I{j}j} 
\]
where $R_j$ is a triangular factor of $S_{\I{j}\I{j}} = R_jR_j^T$.
\end{enumerate}
\end{list}
\end{algdesc}

The algorithm for evaluating $\mathcal R^\mathrm{adj}$ consists 
of the second half of step~2 of Algorithm~\ref{a-hessian} and
of step~3. It also readily follows by taking the adjoint of 
the calculations in 
Algorithm~\ref{a-hessian-factor}.

\begin{algdesc}{Evaluation of adjoint Hessian factor}
\begin{list}{}{}
\item[\textbf{Input.}]
A matrix $W\in\symm^n_V$, the Cholesky factors $L$, $D$ of a 
positive definite $X\in\SV$, and the projected inverse $S = \proj(X^{-1})$.
\item[\textbf{Output.}]
The matrix $T = \mathcal R^\mathrm{adj}(W)$ 
where $\mathcal H = \mathcal R^\mathrm{adj} \circ \mathcal R$ is the 
Hessian of $f$ at $X$.
\item[\textbf{Algorithm.}] \mbox{}
\begin{enumerate}
\item For all $j \in \{1,2,\ldots,n\}$, compute column $j$ of $M$ via
\[
 M_{jj} = \frac{W_{jj}}{D_{jj}}, \qquad 
 M_{\I{j}j} = \frac{1}{\sqrt{D_{jj}}} R_j W_{\I{j}j}. 
\]

\item Iterate over $j\in\{1,2,\ldots,n\}$ in reverse topological 
order.  For each $j$, calculate the $j$th column of $T$ from
\[
 \left[\begin{array}{cc}
  T_{jj} & T_{\I{j}j}^T \\ T_{\I{j}j} & V_j' \end{array}\right] =
  \left[\begin{array}{cc} 1 & -L_{\I{j}j}^T  \\ 0 & I \end{array}\right]
  \left[\begin{array}{cc}
  M_{jj}  &  M_{\I{j}j}^T \\ M_{\I{j}j} & V_j' \end{array}\right] 
  \left[\begin{array}{cc} 1 & 0 \\ -L_{\I{j}j} & I \end{array}\right]
 \]
 and the update matrices $V_i'$ for the children of vertex $j$ via
 \[
  V_i' = E_{\Ip{j}\I{i}}^T \left[\begin{array}{cc}  
  T_{jj} & T_{\I{j}j}^T \\ 
  T_{\I{j}j} & V_j' \end{array}\right] E_{\Ip{j}\I{i}}, \quad
  i \in \ch(j).
 \]
\end{enumerate}
\end{list}
\end{algdesc}

\subsection{Sparse arguments}
In many applications, such as interior-point methods for the
cone programs~(\ref{e-cone-LPs}) mentioned in the introduction, 
the Hessian $\mathcal H = \nabla^2 f(X)$ and its factorization
$\mathcal H = \mathcal R^\mathrm{adj} \circ \mathcal R$ are needed to 
compute coefficients
\[
 H_{ij} = A_i \sdot  (\nabla f(X)[A_j]) 
        = \mathcal R(A_i) \sdot \mathcal R(A_j)
\]
for $m$ matrices $A_1, \ldots, A_m \in\SV$.
The matrices $A_i$ are often very sparse relative to the
sparsity pattern $V$. 
In this section, we examine the implications of sparsity in the matrix 
$Y$ on the computation of $W = \mathcal R(Y)$ using
Algorithm~\ref{a-hessian-factor}.

From step~1 in Algorithm~\ref{a-hessian-factor} we see that if $Y_{jj}
\neq 0$, then $K_{\I{j}j}$ and $U_j'$ are nonzero and dense.
In step~1 these nonzeros are then further propagated via the recursion in 
topological order to all the columns indexed by ancestors of $j$.  
Therefore, if $Y_{jj} \neq 0$, then $W_{\Ip{k}k} \neq 0$ are dense for all 
ancestors $k$ of $j$.

To see how an off-diagonal nonzero in $Y$ affects the sparsity
pattern of $W$, suppose that $Y_{jj} = 0$, $Y_{ij} \neq 0$ for some
$i \in \I{j}$, and $U_k' = 0$ for all $k \in \ch(j)$. 
Then $W_{pq} \neq 0$ for all ancestors $q$ of $j$
and $p\in \Ip{j} \cap \{i,i+1,\ldots,n\}$.
Hence nonzeros in column $j$ of $Y$ create fill in the
columns that correspond  to ancestors of $j$.
In Algorithm~\ref{a-hessian-factor}, we can therefore
prune the elimination tree at node $k$ if all the descendants of 
node $k$ correspond to columns in $Y$ with no lower triangular 
nonzeros.

\begin{figure}
\centering
\begin{tikzpicture}[scale=0.400000]

      \tikzset{VertexStyleA/.style = {shape = circle, draw, minimum size = 5pt, inner sep = 0 pt, color=black}}

      \tikzset{VertexStyleB/.style = {shape = circle, draw, minimum size = 5pt, inner sep = 0 pt, fill=black}}  

      \tikzset{VertexStyleC/.style = {shape = circle, draw, minimum size = 5pt, inner sep = 0 pt, fill=black!25}}  

      \tikzset{VertexStyleD/.style = {shape = rectangle, draw, minimum size = 5pt, inner sep = 0 pt, fill=black}} 

      \foreach \i/\j in {2/0,4/2,14/2,4/3,14/3,8/4,14/4,15/4,8/5,15/5,7/6,8/6,14/6,8/7,14/7,14/8,15/8,10/9,12/9,13/9,16/9,12/10,13/10,16/10,12/11,13/11,15/11,16/11,13/12,15/12,16/12,15/13,16/13,15/14,16/14,16/15}
     \node[VertexStyleA] at (\j,-\i){};

      \foreach \i/\j in {0/0,5/5,6/6,7/7,9/9,10/10,11/11,12/12,13/13}
     \node[VertexStyleA] at (\j,-\i){};

      \foreach \i/\j in {2/1,3/1,3/2,4/2,14/2,4/3,14/3,8/4,14/4,15/4,14/8,15/8,15/14,16/14,16/15}
     \node[VertexStyleB] at (\j,-\i){};

   \foreach \k/\i in {3/2, 4/3, 5/4, 9/8, 15/14, 16/15, 17/16} 
     \node[font=\small] at (\i,-\i) {$\k$};

      \foreach \i/\j in {1/1}
     \node[VertexStyleD] at (\j,-\i){};

      \draw (-0.5,0.5) rectangle (16.5,-16.5);
 
\end{tikzpicture}
\begin{tikzpicture}[scale=0.400000]

      \tikzset{VertexStyleA/.style = {shape = circle, draw, minimum size = 5pt, inner sep = 0 pt, color=black}}

      \tikzset{VertexStyleB/.style = {shape = circle, draw, minimum size = 5pt, inner sep = 0 pt, fill=black}}  

      \tikzset{VertexStyleC/.style = {shape = circle, draw, minimum size = 5pt, inner sep = 0 pt, fill=black!25}}  

      \tikzset{VertexStyleD/.style = {shape = rectangle, draw, minimum size = 5pt, inner sep = 0 pt, fill=black}}  

      \foreach \i/\j in {2/0,2/1,3/1,3/2,4/2,14/2,4/3,14/3,8/4,14/4,15/4,8/5,15/5,7/6,8/6,14/6,8/7,14/7,14/8,15/8,10/9,12/9,13/9,16/9,12/10,13/10,16/10,12/11,13/11,15/11,16/11,13/12,15/12,16/12,15/13,16/13,15/14,16/14,16/15}
     \node[VertexStyleA] at (\j,-\i){};

      \foreach \i/\j in {0/0,1/1,2/2,5/5,6/6,7/7,9/9,10/10,11/11,12/12,13/13}
     \node[VertexStyleA] at (\j,-\i){};

      \foreach \i/\j in {3/2,4/3,14/3,8/4,14/4,15/4,14/8,15/8,15/14,16/14,16/15}
     \node[VertexStyleB] at (\j,-\i){};

   \foreach \k/\i in {4/3, 5/4, 9/8, 15/14, 16/15, 17/16} 
     \node[font=\small] at (\i,-\i) {$\k$};

      \foreach \i/\j in {4/2,14/2}
     \node[VertexStyleC] at (\j,-\i){};

      \foreach \i/\j in {3/1}
     \node[VertexStyleD] at (\j,-\i){};

      \draw (-0.5,0.5) rectangle (16.5,-16.5);
 
\end{tikzpicture}\caption{Sparsity pattern of the lower triangular part of
  $W = \mathcal R(Y)$ when (a) $Y$ has a single diagonal nonzero in 
  position (2,2), marked with a solid square, and (b) when $Y$ has an 
  off-diagonal nonzero in position (4,2). The solid entries mark the 
  nonzero elements in $W$.  The gray markers in (b) correspond to nonzeros 
  in $W$ that are introduced by the scaling in step 2 of
  Algorithm~\ref{a-hessian-factor}, \ie, these entries are not present in 
  the intermediate variable $K$.}
  \label{fig-spy-W}
\end{figure}
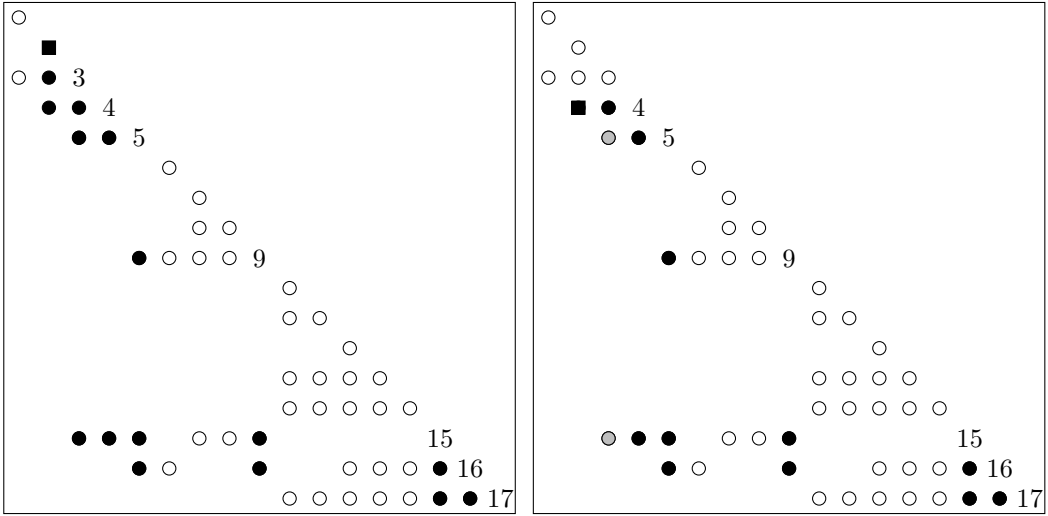

Figure~\ref{fig-spy-W} shows examples of the sparsity pattern of $W$ when 
$Y$ has a diagonal and an off-diagonal nonzero, respectively, for
the pattern $V$ in Figure~\ref{f-etree1}.

\paragraph{Numerical results}
To evaluate the benefits of exploiting additional sparsity in the
argument, we have implemented and tested a version of
Algorithm~\ref{a-hessian-factor} that exploits sparsity in $Y$ 
relative to $V$.
In the experiment, we use as test data a set of
randomly generated problems with sparsity patterns from the
University of Florida Sparse Matrix Collection. For each sparsity
pattern, we generate ten sparse arguments with just two
lower-triangular nonzero entries in random
positions. Table~\ref{tbl:ex_sparse_hessian_factor} summarizes the
average time required to compute $W = \mathcal{R}(Y)$ using
Algorithm~\ref{a-hessian-factor} with and without the techniques
described in this section. For these very sparse arguments, pruning
the elimination tree results in average speedups in the range
4--6, but in general the speedup depends both on the number of
nonzeros and on the position of the nonzeros.

Finally, we remark that additional computational savings can be made
in step 2 of Algorithm~\ref{a-hessian-factor} when $\mathcal{R}(Y)$ is
needed for several arguments $Y$. Specifically, the triangular factors
$R_j$ of $S_{\I{j}\I{j}} = R_jR_j^T$ need only be computed once.

\begin{table}
  \centering
  \begin{tabular}{lrd{2}d{2}d{2}}
    \hline
    Sparsity pattern & \multicolumn{1}{c}{$n$} & \multicolumn{1}{c}{Dense} &  \multicolumn{1}{c}{Sparse} & \multicolumn{1}{c}{Ratio} \\
  &     & \multicolumn{1}{c}{(seconds)}&  \multicolumn{1}{c}{(seconds)} &          \\
    \hline
 HB/plat1919        &  1919 &   0.17 &    0.05 &   3.7  \\   HB/bcsstk13        &  2003 &   1.40 &    0.36 &   3.9  \\   HB/lshp3025        &  3025 &   0.21 &    0.05 &   4.6  \\   Boing/nasa4704     &  4704 &   1.05 &    0.27 &   4.0  \\   TKK/g3rmt3m3       &  5357 &   1.43 &    0.35 &   4.1  \\   Schenk\_IBMNA/c-36 &  7479 &   0.37 &    0.06 &   6.1  \\   Wang/swang1        & 10800 &   7.05 &    1.72 &   4.1  \\   ACUSIM/Pres\_Poisson & 14822 & 17.65 &   4.35 &   4.1  \\   GHS\_psdef/wathen100 & 30401 &  6.05 &   1.38 &   4.4  \\      \hline
  \end{tabular}
  \caption{Average computational times (10 trials) for 
    Algorithm~\ref{a-hessian-factor} with and without the 
    technique for sparse arguments (columns `sparse' and `dense', 
    respectively) when applied to sparse arguments. The sparse 
    arguments are randomly generated with two lower-triangular 
    nonzero entries in random positions. }
  \label{tbl:ex_sparse_hessian_factor}
\end{table}

\section{Cliques and clique trees} \label{s-cliques}
It is known that the performance of sparse Cholesky factorization
algorithms on modern computers can be improved by combining groups of 
vertices into \emph{supernodes} and applying block elimination to the 
corresponding columns.
Several definitions of supernodes exist in the literature.  
In this paper, we define a supernode as a maximal group of 
columns of $L$ (sorted, but not necessarily contiguous) that share the 
same nonzero structure.  More specifically, if $N$ is a supernode and
$j = \max N$, then for all $k\in N$,  
\[
 \Ip{k} = (N \cup \I{j}) \cap \{k, k+1, \ldots,n\}.
\]
In the example of Figure~\ref{f-etree1}, the sets
\[
\{1\}, \quad \{2\}, \quad \{3,4\}, \quad \{5,9\}, \quad \{6\}, \quad
\{7,8\}, \quad \{10,11\}, \quad \{12,13,14\}, \quad \{15,16,17\}
\]
form supernodes.  (Another more common definition 
adds the requirement that the indices in $N$ are contiguous
\cite{LNP:93};
this can be achieved from a set of supernodes as defined above by
a simple reordering.)  
In the context of multifrontal factorizations, the grouping into
supernodes has the advantage that only one frontal matrix
is required per supernode.  This reduces the memory and arithmetic
overhead incurred for the assembly of frontal matrices.
Moreover, the block operations allow us to replace matrix-vector 
operations (level-2 BLAS) with more efficient matrix-matrix operations 
(level-3 BLAS) \cite{DDHD:90}.

Supernodes are closely related to cliques in the filled graph and
the barrier algorithms described in \cite{DVR:08,DaV:09a},
which involve iterations on clique trees, can be interpreted as supernodal 
multifrontal algorithms.  We therefore start the discussion with
a review of cliques and clique trees, and their connections with
supernodes.

\subsection{Cliques}
A filled graph is also known as a \emph{chordal} or \emph{triangulated}
graph.
A \emph{clique} is a maximal set of vertices that define a complete 
subgraph of the filled graph.  Equivalently, a clique is a set of 
indices that define a dense lower-triangular principal subblock of $L$.
Every clique $W$ in a filled graph can be expressed as $W = \Ip{i}$, 
where $i = \min W$, the least element in $W$ \cite[proposition 2]{BlP:93}.
This follows from the fact that the index sets $\Ip{i}$ define complete 
subgraphs, as noted in Section~\ref{s-etree}.
Hence, if $i$ is the lowest index in the clique $W$, then 
$W\subseteq \Ip{i}$.  Since $W$ is maximal, we must have $W = \Ip{i}$.
The vertex $i = \min W$ is called the \emph{representative vertex} of the 
clique.
Since there are at most $n$ representative vertices, a filled graph can 
have at most $n$ cliques.

Efficient algorithms for identifying the representative vertices
can be derived from the following criterion:  the cliques are 
exactly the sets $\Ip{i}$ for which there exists no $j<i$ with
$\Ip{i} \subset \Ip{j}$ \cite[proposition 3]{LPP:89}.
This follows from the characterization of cliques in terms of 
representative vertices. 
If $\Ip{i} \subset \Ip{j}$ for some $j<i$, then $\Ip{i}$ is certainly
not a clique, since it is strictly included in another complete
subgraph.
Conversely, if $\Ip{i}$ is not a clique, \ie, $\Ip{i} \subset W$  for
some clique $W$, and $j$ is the representative vertex of $W$, then
$j<i$ and $\Ip{i} \subset W = \Ip{j}$.

In the example in Figure~\ref{f-etree1}, the representative vertices 
are $1$, $2$, $3$, $5$, $6$, $7$, $10$, $12$, $15$.  
The other vertices are not representative  because the corresponding sets 
$\Ip{j}$ are not maximal:
\[
\Ip{4} \subset \Ip{3}, \quad \Ip{8} \subset \Ip{7},
\quad \Ip{9} \subset \Ip{5}, \quad \Ip{11} \subset \Ip{10},
\quad \Ip{13}, \Ip{14}, \Ip{16}, \Ip{17} \subset \Ip{12}.
\]
The representative vertices are easily identified from the elimination 
tree and the monotone degrees of the vertices.  It can be shown that a 
vertex $j$ is a representative vertex if and only if 
\BEQ \label{e-repr-test}
  |\I{j}| > |\I{k}| -1 \quad \forall k\in\ch(j) 
\EEQ
(see \cite{PoS:90}).  To see this, recall from~(\ref{e-etree-ch})
that $\I{k} \subseteq \Ip{j}$ and $|\I{j}| \geq |\I{k}|-1$ if 
$k\in \ch(j)$.
Therefore, if $|\I{j}| = |\I{k}| -1$ for some $k\in\ch(j)$,
then $\Ip{j} = \I{k} \subset \Ip{k}$.  Therefore $j$ is not a 
representative 
vertex because the complete subgraph defined by $\Ip{j}$ is not maximal.
Conversely, suppose $j$ is not representative, \ie, 
$\Ip{j} \subset \Ip{l}$ for some $l<j$.  In particular, 
$L_{jl} \neq 0$ and therefore, $l$ is a descendant of $j$ in the 
elimination tree.
From Theorem~\ref{t-column-dep}
(the set of inequalities~(\ref{e-etree-path})), this implies
that $\Ip{j} \subset \Ip{k}$ for all $k$ in the path from $l$ to $j$. 
In particular, 
\[
 \Ip{j} \subset \Ip{k} \subseteq \Ip{j} \cup \{k\}
\] 
for the child $k$ of $j$ on this path.  
Therefore, $\Ip{j} = \I{k}$ and $|\I{j}| = |\I{k}|-1$.

\subsection{Supernode partitions} \label{s-supernodes}
By comparing the monotone degrees of the vertices in the elimination tree 
and their parents we can partition the vertices $\{1,2,\ldots,n\}$ in sets 
$\new{j}  = \{i_1, i_2, \ldots, i_r\}$ where $i_1=j$ is a representative 
clique vertex. The vertices $i_1$, $i_2$, \ldots, $i_r$ form a path 
from $j$ to an ancestor $i_r$ of $j$ in the elimination tree, and 
\[
 |\I{i_1}| \; = \; |\I{i_2}| + 1 \; = \; |\I{i_3}| + 2 \; = \; \cdots 
 \; = \; |\I{i_r}| + r -1,
\]
or, equivalently,
\BEQ \label{e-new2}
 \I{i_1} \; = \; \Ip{i_2} 
 \;  = \; \{i_2\} \cup \Ip{i_3}  \; =  
 \;  \cdots \; = \; \{i_2, i_3, \ldots, i_{r-1} \} \cup \Ip{i_r}.
\EEQ
The sets $\new{j}$ are supernodes (in the definition given at the
beginning of this section).
In general, several such partitions exist.
Two possible partitions for the elimination tree
in Figure~\ref{f-etree1} are shown in Figure~\ref{f-etree} 
and listed in  Table~\ref{t-example}.
The representative vertices are shown as rectangles, and the sets 
$\new{j}$ are the vertices on the paths shown with heavy lines.   
\begin{figure}

\hspace*{\fill}
\begin{tikzpicture}[scale=10.00]

      \tikzset{VertexStyleA/.style = {shape = circle, draw, 
       minimum size = 1, inner sep = 2, fill=none}}
   \tikzset{VertexStyleB/.style = {shape = rectangle, draw, 
       minimum size = 1, inner sep = 4, fill=none}}

      \foreach \l/\k/\d/\x/\y in {
       1/1/1/0.13/0.0,
       2/2/2/0.27/0.0,
       3/3/3/0.2/0.1,
       5/5/3/0.2/0.3,
       6/6/2/0.35/0.3,
       7/7/3/0.5/0.2,
       10/10/4/0.82/0.2,
       12/12/4/0.68/0.3,
       15/15/2/0.35/0.5}
       \node[VertexStyleB, font=\small](\k) 
           [label = right: \footnotesize $\d$] 
           at (\x,\y){\l};

      \foreach \l/\k/\d/\x/\y in {
       4/4/2/0.2/0.2,
       8/8/2/0.5/0.3,
       9/9/2/0.35/0.4,
       11/11/3/0.82/0.3,
       13/13/3/0.75/0.4,
       14/14/2/0.75/0.5,
       16/16/1/0.55/0.6,
       17/17/0/0.55/0.7}
       \node[VertexStyleA, font=\small](\k) 
           [label = right: \footnotesize $\d$] 
           at (\x,\y){\l};

      \foreach \i/\j in {
       1/3,
       2/3,
       4/5,
       6/9,
       8/9,
       9/15,
       14/16,
       11/13}
       \draw[thin, dashed] (\i)--(\j);

   \foreach \i/\j in {
       3/4,
       5/9,
       7/8,
       10/11,
       12/13,
       13/14,
       15/16,
       16/17}
       \draw[ultra thick] (\i)--(\j);

\end{tikzpicture}
\hspace*{\fill}
\hspace*{\fill}
\begin{tikzpicture}[scale=10.00]

      \tikzset{VertexStyleA/.style = {shape = circle, draw, 
       minimum size = 1, inner sep = 2, fill=none}}
   \tikzset{VertexStyleB/.style = {shape = rectangle, draw, 
       minimum size = 1, inner sep = 4, fill=none}}

      \foreach \l/\k/\d/\x/\y in {
       1/1/1/0.13/0.0,
       2/2/2/0.27/0.0,
       3/3/3/0.2/0.1,
       5/5/3/0.2/0.3,
       6/6/2/0.35/0.3,
       7/7/3/0.5/0.2,
       10/10/4/0.82/0.2,
       12/12/4/0.68/0.3,
       15/15/2/0.35/0.5}
       \node[VertexStyleB, font=\small](\k) 
           [label = right: \footnotesize $\d$] 
           at (\x,\y){\l};

      \foreach \l/\k/\d/\x/\y in {
       4/4/2/0.2/0.2,
       8/8/2/0.5/0.3,
       9/9/2/0.35/0.4,
       11/11/3/0.82/0.3,
       13/13/3/0.75/0.4,
       14/14/2/0.75/0.5,
       16/16/1/0.55/0.6,
       17/17/0/0.55/0.7}
       \node[VertexStyleA, font=\small](\k) 
          [label = right: \footnotesize $\d$] 
           at (\x,\y){\l};

      \foreach \i/\j in {
       1/3,
       2/3,
       4/5,
       6/9,
       8/9,
       9/15,
       15/16,
       11/13}
       \draw[thin, dashed] (\i)--(\j);

   \foreach \i/\j in {
       3/4,
       5/9,
       7/8,
       10/11,
       12/13,
       13/14,
       14/16,
       16/17}
       \draw[ultra thick] (\i)--(\j);

\end{tikzpicture}
\hspace*{\fill}
\caption{Two supernode partitions of the elimination tree in
Figure~\ref{f-etree1}.  The representative vertices are
shown as rectangles.  The vertices joined by the paths shown with
heavy lines form the sets $\new{k}$, where $k$ is the representative 
vertex on the path.   The sets $\new{k}$ are enumerated in 
Table~\ref{t-example}.  \label{f-etree}} 
\end{figure}\begin{table}
\begin{center}
\begin{tabular}{c|l|l} 
$k$ & \multicolumn{1}{c|}{$\new{k}$}  
& \multicolumn{1}{c}{$\anc{k}$} \\\hline
  1 & $\{1\}$           & $\{3\}$        \\
  2 & $\{2\}$           & $\{3,4\}$      \\
  3 & $\{3,4\}$         & $\{5,15\}$     \\
  5 & $\{5, 9\}$        & $\{15,16\}$    \\
  6 & $\{6\}$           & $\{9,16\}$     \\
  7 & $\{7,8\}$         & $\{9,15\}$     \\
 10 & $\{10,11\}$       & $\{13,14,17\}$ \\
 12 & $\{12,13,14\}$    & $\{16,17\}$    \\
 15 & $\{15,16,17\}$    & $\{\}$          
\end{tabular} 
\qquad\qquad
\begin{tabular}{c|l|l} 
$k$ & \multicolumn{1}{c|}{$\new{k}$}  
& \multicolumn{1}{c}{$\anc{k}$} \\\hline
  1 & $\{1\}$                 & $\{3\}$       \\
  2 & $\{2\}$                 & $\{3,4\}$     \\
  3 & $\{3,4\}$               & $\{5,15\}$    \\
  5 & $\{5, 9\}$              & $\{15,16\}$    \\
  6 & $\{6\}$                 & $\{9,16\}$     \\
  7 & $\{7,8\}$               & $\{9, 15\}$    \\
 10 & $\{10,11\}$             & $\{13,14,17\}$ \\
 12 & $\{12,13,14,15,16,17\}$ & $\{\}$         \\
 15 & $\{15\}$                & $\{16,17\}$ 
\end{tabular}
\end{center}
\caption{The two supernode partitions defined in Figure~\ref{f-etree}.
The first columns of the tables are the representative vertices.
Each clique $\Ip{k}$  is partitioned in two sets as 
$\Ip{k} = \new{k} \cup \anc{k}$.  The sets $\new{k}$ for a partition
 of the vertices $\{1,2,\ldots,n\}$} \label{t-example}
\end{table}
We note two important properties of the sets $\new{j}$:
\BIT
\item The set $\new{j}$ is a subset of the clique represented by 
 vertex $j$: $\new{j}\subseteq \Ip{j}$.  This can be seen 
from~(\ref{e-new2}) which implies $\Ip{i} \subset \Ip{j}$ if 
 $i\in\new{j}$ and $i\neq j$.

\item Define $\anc{j} = \Ip{j} \setminus \new{j}$.
Then we have $i< \min\anc{j}$ for all $i\in\new{j}$.
To see this, first note that if $k\in\Ip{j}$, then $L_{kj} \neq 0$ 
 and 
therefore, $k$ is an ancestor of $j$ in the elimination tree 
(Theorem~\ref{t-ancestor}).
If also $k\leq i$ for some $i\in\new{j}$, then $k$ is on the path from 
$j$ to $i$ in the elimination tree.  
However, by definition of $\new{j}$, this means that $k\in\new{j}$.
\EIT
This result means that $\new{j} \cup \anc{j}$ is an ordered
partition of the clique $\Ip{j}$, \ie, the elements of $\new{j}$
have a lower index than the elements of $\anc{j}$.
(In \cite{LPP:89,PoS:90}, the sets 
$\new{j}$ and $\anc{j}$ are referred to as the \emph{new set}
$\mbox{new}(K)$ and the \emph{ancestor set} 
$\mbox{anc}(K)$, respectively, where $K$ is the clique
$K = \Ip{j}$.)

If $\anc{j}$ is nonempty, we refer to the vertex 
$k= \min\anc{j}$ as the \emph{first ancestor} of 
the clique $\Ip{j}$.  
The first ancestor can be identified from the elimination tree
as the parent of the vertex $i = \max\new{j}$.
This follows from~(\ref{e-new2}) with $i=i_r$ and the fact that the
parent of vertex $i$ is the first element in $\Ip{i}$ greater than $i$.
Note that while the first ancestor can be determined from
the elimination tree and the sets $\new{j}$, the rest of the 
sets $\anc{j}$ cannot be derived from the elimination tree but require
knowledge of the clique $\Ip{j}$.

The sets $\new{k}$ and $\anc{k}$ for the two partitions in the 
example are listed in Table~\ref{t-example}.

\subsection{Clique trees} \label{s-clique-trees}
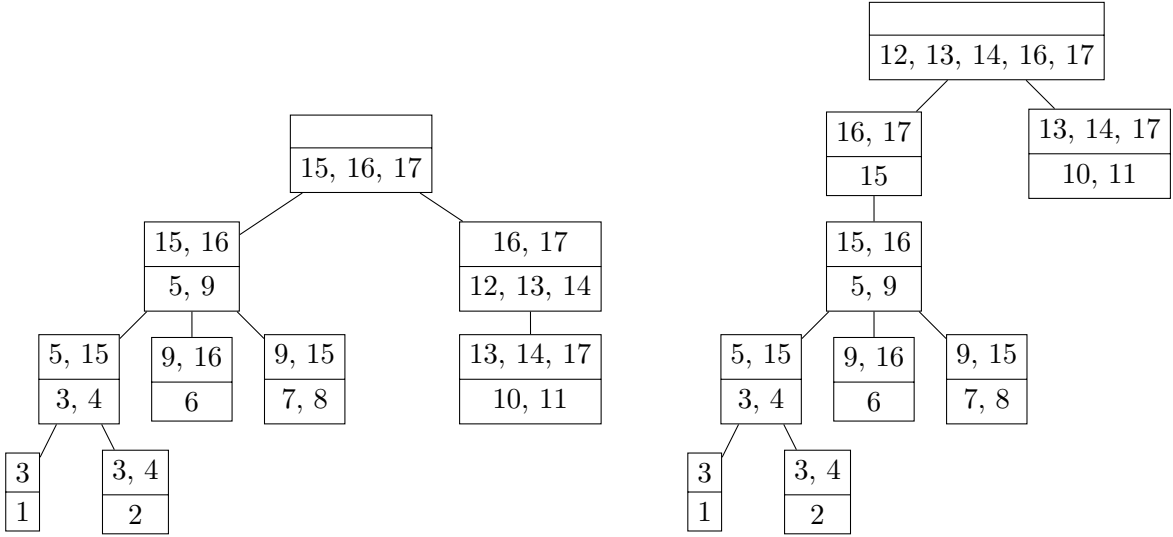
\begin{figure}

\hspace*{\fill}
\begin{tikzpicture}[scale=1.5,
    clique/.style = {rectangle split, rectangle split parts=2, draw},
    clique2/.style = {shape = circle split, draw}]
 \node[clique](1) at (1.5,-1){3 \nodepart{second} 1 };
 \node[clique](2) at (2.5,-1){3, 4 \nodepart{second} 2};
 \node[clique](3) at (2,0){5, 15 \nodepart{second} 3, 4};
 \node[clique](6) at (3,0){9, 16 \nodepart{second} 6};
 \node[clique](7) at (4,0){9, 15 \nodepart{second} 7, 8};
 \node[clique](5) at (3,1){15, 16 \nodepart{second} 5, 9};
 \node[clique](10) at (6,0){13, 14, 17 \nodepart{second} 10, 11};
 \node[clique](12) at (6,1){16, 17 \nodepart{second} 12, 13, 14};
 \node[clique](15) at (4.5,2){ \nodepart{second} 15, 16, 17};

 \draw (1) -- (3);
 \draw (2) -- (3);
 \draw (3) -- (5);
 \draw (6) -- (5);
 \draw (7) -- (5);
 \draw (5) -- (15);
 \draw (10) -- (12);
 \draw (12) -- (15);
\end{tikzpicture}
\hspace*{\fill}
\hspace*{\fill}
\begin{tikzpicture}[scale=1.5,
    clique/.style = {rectangle split, rectangle split parts=2, draw},
    clique2/.style = {shape = circle split, draw}]
 \node[clique](1) at (1.5,-1){3 \nodepart{second} 1 };
 \node[clique](2) at (2.5,-1){3, 4 \nodepart{second} 2};
 \node[clique](3) at (2,0){5, 15 \nodepart{second} 3, 4};
 \node[clique](6) at (3,0){9, 16 \nodepart{second} 6};
 \node[clique](7) at (4,0){9, 15 \nodepart{second} 7, 8};
 \node[clique](5) at (3,1){15, 16 \nodepart{second} 5, 9};
 \node[clique](15) at (3,2){16, 17 \nodepart{second} 15};
 \node[clique](10) at (5,2){13, 14, 17 \nodepart{second} 10, 11};
 \node[clique](12) at (4,3){ \nodepart{second} 12, 13, 14, 16, 17};

 \draw (1) -- (3);
 \draw (2) -- (3);
 \draw (3) -- (5);
 \draw (6) -- (5);
 \draw (7) -- (5);
 \draw (5) -- (15);
 \draw (15) -- (12);
 \draw (10) -- (12);

\end{tikzpicture}
\hspace*{\fill}

\caption{The two clique trees corresponding to the partitions
in Figure~\ref{f-etree}.  
The first index in the bottom row at each node is the representative
vertex $j$ of the clique $\Ip{j}$.  
Each clique $\Ip{j}$ is partitioned in two sets 
$\Ip{j} = \new{j} \cup \anc{j}$.
The indices of the bottom row at each node form the sets $\new{j}$; 
the indices in the top row form the set $\anc{j}$.
The parent of clique $\Ip{j}$ is the clique $\Ip{k}$ that includes the 
first ancestor of clique $\Ip{j}$ in its $\new{k}$ set.
\label{f-etree2}}
\end{figure}

We can associate with the vertex partitioning in sets $\new{j}$ 
a tree with the cliques $\Ip{j}$ as its nodes.
The root of the clique tree is the clique represented 
by the vertex $j$ for which $n\in\new{j}$.
The parent of the clique $\Ip{j}$ is the clique which has as its
representative the vertex $k$ for which $\min\anc{j} \in \new{k}$. 
We will use the notation $k=\pa{j}$ to denote that $\Ip{k}$
is the parent of $\Ip{j}$ in the clique tree.
Figure~\ref{f-etree2} shows the  clique trees defined by the 
partitions $\new{j}$ in Figure~\ref{f-etree}.
The clique tree satisfies the following key properties 
\cite[p.186]{PoS:90} \cite{LPP:89}:
\BIT
\item $\anc{j} \subset \Ip{\pa{j}}$.

Indeed, let $i= \min\anc{j}$ be the first ancestor of clique
$\Ip{j}$.   By definition of the clique tree, $i\in\new{\pa{j}}$.
From the definition in~(\ref{e-new2}) this implies that
$\Ip{i} \subset \Ip{\pa{j}}$.
Since $\anc{j}$ defines a complete subgraph of vertices with indices 
greater than or equal to $i$, we have $\anc{j} \subseteq \Ip{i}$.
Therefore $\anc{j} \subset \Ip{\pa{j}}$.

\item An element of $\new{j}$ is in the clique $\Ip{k}$ only if $\Ip{k}$ 
is a descendant of $\Ip{j}$ in the clique tree.

We can show this by contradiction. 
Suppose $i\in\new{j}$ belongs to $\Ip{k}$ and $\Ip{k}$ is not
a descendant of $\Ip{j}$ in the clique tree.
The sets $\new{l}$ form a partition of $\{1,2,\ldots,n\}$, so if
$i\in\new{j}$ and $i\in\Ip{k}$ for $k\neq j$, then $i\in\anc{k}$.  
By the previous property, this implies $i\in\Ip{\pa{k}}$.
We have $\pa{k} \neq j$ because $\Ip{k}$ is not a descendant of $\Ip{j}$.
Therefore $i\in\anc{\pa{k}}$ and, again from the previous property,
$i\in\Ip{\pa{\pa{k}}}$ and  $i\in\anc{\pa{\pa{k}}}$.
Continuing this process recursively, we eventually arrive at the 
conclusion that $i$ belongs to $\anc{r}$ where $\Ip{r}$ is the
root of the clique tree.  However, this is impossible
because $\Ip{r} = \new{r}$ and $\anc{r} = \emptyset$.
 
\item If an element of $\new{j}$ is in $\anc{k}$, then it belongs to
all the cliques on the path between $\Ip{j}$ and $\Ip{k}$
in the clique tree.

This follows by combining the first two properties.  From the
second property, if $i\in\new{j}$ and $i\in\anc{k}$, then 
$\Ip{k}$ is a descendant of $\Ip{j}$.
Assume there are cliques $\Ip{s}$, $\Ip{\pa{s}}$ on the path between 
$\Ip{k}$ and $\Ip{j}$ with the property that $i\in\Ip{s}$ 
and $i\not\in\Ip{\pa{s}}$.  From the first property, this implies
$i\in\new{s}$.  But this contradicts $i\in\new{j}$, unless $j=s$,
because $\new{j} \cup \new{s} = \emptyset$ if $j\neq s$.

\EIT
Taken together, these three properties state that the cliques that
contain a vertex $i$ form a subtree in the clique tree.
The root of the subtree is the unique clique $\Ip{k}$ for 
which $i\in\new{k}$.   This is known as the \emph{induced subtree
property} of clique trees \cite{BlP:93}.

\subsection{Clique tree algorithm}
To summarize the results of this section, we state  a simple
algorithm that identifies the representative vertices of the cliques,
generates a partition into sets $\new{k}$, identifies the first 
ancestors $\min\anc{k}$ of the cliques, and determines the parent 
structure of the clique tree.  The algorithm is due to
Pothen and Sun \cite[p.185]{PoS:90}.

\begin{algdesc}{Clique tree algorithm} 
\begin{list}{}{}
\item[\textbf{Input.}] An elimination tree and the monotone degree
$|I_k|$, $k=1,\ldots,n$. 
\item[\textbf{Output.}] The representative vertices, the partition
 in supernodes $\new{j}$,  the first ancestor of each clique, 
 and the parent structure of a clique tree.
\newpage
\item[\textbf{Algorithm.}]
For $i=1,\ldots, n$:
\begin{enumerate}
\item If $|\I{i}| > |\I{j}|-1$ for all $j \in\ch(i)$, then $i$ is 
 a representative vertex.  Set $\new{i} = \{i\}$ and $k=i$.
 Otherwise, choose a vertex $j\in\ch(i)$ with $|\I{i}| = |\I{j}|-1$, 
 determine the representative vertex $k$ for which $j\in\new{k}$,
 and add $i$ to $\new{k}$.

\item For each $j\in \ch(i)$, if $j \in \new{l}$ and $l\neq k$,  
 set $\pa{l} := k$  and $\min\anc{l} := i$.
\end{enumerate}
\end{list}
\end{algdesc}

Note that in step~1, there may be several choices for 
the child vertex $j$, and these choices lead to different vertex 
partitions and different clique trees.
The vertex $i=16$ in Figure~\ref{f-etree}, for example, has two
children $j$ that both satisfy $|\I{i}| = |\I{j}| - 1$.
These two choices lead to the different vertex partitions in
Figure~\ref{f-etree} and the two clique trees in Figure~\ref{f-etree2}.

\section{Supernodal multifrontal algorithms} \label{s-supernodal}
We assume there are $l$ cliques, with representative nodes
$i_1$, \ldots, $i_l$
and that the sets $\new{i}$ are contiguous.
The supernodal algorithms are block versions of the multifrontal 
algorithms in which the scalar diagonal elements $X_{jj}$ are replaced with
dense principal blocks $X_{\new{j}\new{j}}$  and the subcolumns
$X_{\I{j}j}$ with dense submatrices $X_{\anc{j}\new{j}}$.

For a clique $\Ip{k}$ we denote by $\ch(\Ip{k})$ the set of 
child cliques of $\Ip{k}$ in the clique tree.
This is not to be confused with $\ch(k)$ (with a vertex $k$ as argument),
which refers to the children of the vertex  $k$ in the elimination tree.

In this section we start with a supernodal version of the 
Cholesky factorization algorithm.
We then give similar extensions of the primal and dual gradient 
evaluation algorithms.
For the sake of brevity, we will omit the extensions of the 
other algorithms in Sections~\ref{s-chol}--\ref{s-hess}, which
follow the same pattern.

\subsection{Cholesky factorization}
In the supernodal Cholesky factorization we factor $X$ as $X = LDL^T$ 
with $D$ block-diagonal and $L$ unit lower triangular.
The matrix $D$ has $l$ dense diagonal blocks $D_{\new{i}\new{i}}$
for $i \in \{i_1,\ldots, i_l\}$.
Corresponding with each clique, $L$ has a diagonal block 
$L_{\new{i}\new{i}} = I$ and a dense submatrix 
$L_{\anc{i}\new{i}}$.  The rest of the block-column indexed by
$\new{i}$ is zero.

\begin{algdesc}{Cholesky factorization} \label{a-sn-chol}
\begin{list}{}{}
\item[\textbf{Input.}]  A positive definite matrix $X\in \SV$ and
 a clique tree for the sparsity pattern $V$.
\item[\textbf{Output.}] The factors $L$, $D$ in the Cholesky factorization
  $X = LDL^T$.
\item[\textbf{Algorithm.}]
Iterate over $j\in\{i_1, i_2, \ldots, i_l\}$ using a topological order 
of the clique tree.   For each $j$, form the frontal matrix 
\[
 F_j = \left[\begin{array}{cc} F_{11} & F_{21}^T \\ F_{21} & F_{22} 
 \end{array}\right] = 
 \left[\begin{array}{cc}
  X_{\new{j}\new{j}} & X_{\anc{j}\new{j}}^T \\ X_{\anc{j}\new{j}} & 0 
 \end{array}\right]
 + \sum_{\Ip{i} \in \ch(\Ip{j})} 
  E_{\Ip{j}\anc{i}} U_i E_{\Ip{j}\anc{i}}^T
\]
and calculate $D_{\new{j}\new{j}}$, $L_{\anc{j}\new{j}}$, and the
update matrix $U_j$ from
\[
 D_{\new{j}\new{j}} = F_{11}, \qquad
 L_{\anc{j}\new{j}} = F_{21} D_{\new{j}\new{j}}^{-1}, \qquad
 U_j = F_{22} - L_{\anc{j}\new{j}} D_{\new{j}\new{j}} L_{\anc{j}\new{j}}^T.
\]
\end{list}
\end{algdesc}
As can be seen, only one frontal matrix is assembled per clique,
a major advantage compared to Algorithm~\ref{a-chol}.
Moreover, the main computation is the level-3 BLAS operation 
in the computation of $U_j$.

\subsection{Gradients}
The supernodal counterpart of Algorithms~\ref{a-projinv} for computing
the primal gradient or projected inverse is as follows.
In this algorithm, $V_j$ is a dense `update matrix' defined
as $V_j = S_{\anc{j}\anc{j}}$.

\begin{algdesc}{Projected inverse} \label{a-sn-projinv}
\begin{list}{}{}
\item[\textbf{Input.}] The Cholesky factors $L$, $D$ of a positive 
  definite matrix $X\in\SV$.
\item[\textbf{Output.}] The projected inverse $S = \proj(X^{-1}) 
 = -\nabla f(X)$.
\item[\textbf{Algorithm.}]  Iterate over $j\in \{i_1, \ldots, i_l\}$
 using a reverse topological order of the clique tree.  
 For each $j$, calculate $S_{\new{j}\new{j}}$ and $S_{\anc{j}\new{j}}$ 
 from
\BEQ \label{e-sn-gradient-1}
 S_{\anc{j}\new{j}} = -V_j L_{\anc{j}\new{j}}, \qquad 
 S_{\new{j}\new{j}} = D_{\new{j}\new{j}}^{-1} - 
    S_{\anc{j}\new{j}}^TL_{\anc{j}\new{j}} 
\EEQ
and compute the update matrices  
\BEQ \label{e-sn-gradient-2}
  V_i = E_{\Ip{j}\anc{i}}^T \left[\begin{array}{cc} 
 S_{\new{j}\new{j}} & S_{\anc{j}\new{j}}^T \\ 
 S_{\anc{j}\new{j}} & V_j \end{array}\right] 
  E_{\Ip{j}\anc{i}}, \quad
 \Ip{i}\in \ch(\Ip{j}). 
\EEQ
\end{list}
\end{algdesc}
The main calculation is the matrix-matrix product 
in~(\ref{e-sn-gradient-1}) which replaces the matrix-vector 
product~(\ref{e-mf-gradient-1}) in the multifrontal algorithm.

The extension of Algorithm~\ref{a-completion} for computing the dual
gradient or the maximum determinant positive definite completion
is as follows.
\begin{algdesc}{Matrix completion}  \label{a-sn-completion}
\begin{list}{}{}
\item[\textbf{Input.}] A matrix $S\in\SV$ that has a positive definite
 completion.
\item[\textbf{Output.}] The Cholesky factors $L$, $D$ of
 $X= -\nabla f_*(S)$, \ie, of the positive definite matrix $X\in\SV$ that
 satisfies $\proj(X^{-1}) = S$.
\item[\textbf{Algorithm.}]
Iterate over $j\in\{i_1,  \ldots, i_l\}$ using a reverse topological
  order of the clique tree.  For each $j$, compute
  $D_{\new{j}\new{j}}$ and $L_{\anc{j}\new{j}}$ from
\BEQ \label{e-sn-compl-1}
 L_{\anc{j}\new{j}} = -V_j^{-1} S_{\anc{j}\new{j}}, \qquad 
 D_{\new{j}\new{j}} = (S_{\new{j}\new{j}} + 
  S_{\anc{j}\new{j}}^T L_{\anc{j}\new{j}})^{-1},
\EEQ
and compute the update matrices
\BEQ \label{e-sn-compl-2}
  V_i = E_{\Ip{j}\anc{i}}^T \left[\begin{array}{cc} 
  S_{\new{j}\new{j}} & S_{\anc{j}\new{j}}^T \\ 
  S_{\anc{j}\new{j}} & V_j \end{array}\right] 
  E_{\Ip{j}\anc{i}}, \quad \Ip{i}\in \ch(\Ip{j}).
\EEQ
\end{list}
\end{algdesc}

As for the multifrontal completion algorithm with factored update
matrices (Algorithm~\ref{a-completion-factored}), 
this algorithm can be improved by propagating a factorization of $V_j$ 
and using~(\ref{e-sn-compl-2}) to compute the factors of $V_i$ from the 
factors of $V_j$.
We mentioned in section~\ref{s-supernodes} that for every clique,
the vertices in $\new{i}$ precede those in $\anc{i}$.  As a consequence,
the matrix $E_{\Ip{j}\anc{i}}$ can be partitioned as
\[ E_{\Ip{j}\anc{i}} =
\left[ \begin{array}{cc}
  E_{\new{j}, \anc{i} \cap \new{j}} & 0 \\
  0 & E_{\anc{j}, \anc{i} \setminus \new{j} }
\end{array} \right].
\]
Using this property in~(\ref{e-sn-compl-2}) we get
\[ V_i =  E_{\Ip{j}\anc{i}}^T \left[\begin{array}{cc} 
  S_{\new{j}\new{j}} & S_{\anc{j}\new{j}}^T \\ 
  S_{\anc{j}\new{j}} & V_j \end{array}\right] 
  E_{\Ip{j}\anc{i}} = \left[ \begin{array}{cc} A & B^T \\ B & CC^T\end{array} \right]
\]
where $A$ is a principal submatrix of $S_{\new{j}\new{j}}$ of order
$|\anc{i}\cap \new{j}|$, and $C$ consists of $| \anc{i} \setminus
\new{j} |$ rows of the upper triangular factor $R_j$ of $V_j =
R_jR_j^T$. By reducing $C$ to square triangular form $C =
RQ^T$ using a series of Householder transformations, 
and a factorization $A - B^TR^{-T}R^{-1}B = \tilde{R}\tilde{R}^T$ 
we obtain the factorization $V_i = R_iR_i^T$ as
\[
  R_i = \left[
    \begin{array}{cc}
      \tilde{R} & (R^{-1}B)^T \\
      0 & R
    \end{array} \right].
\]

\subsection{Numerical results}
\begin{figure}
  \centering
  \subfloat[]{\label{fig-ufsmc_sn_chol-vs-sn_parinv}
\psfrag{x00}[c][c]{${10^{-3}}$}
\psfrag{x01}[c][c]{${10^{-2}}$}
\psfrag{x02}[c][c]{${10^{-1}}$}
\psfrag{x03}[c][c]{${10^{0}}$}
\psfrag{x04}[c][c]{${10^{1}}$}
\psfrag{x05}[c][c]{${10^{2}}$}
\psfrag{y00}[r][r]{${10^{-3}}$}
\psfrag{y01}[r][r]{${10^{-2}}$}
\psfrag{y02}[r][r]{${10^{-1}}$}
\psfrag{y03}[r][r]{${10^{0}}$}
\psfrag{y04}[r][r]{${10^{1}}$}
\psfrag{y05}[r][r]{${10^{2}}$}
\psfrag{xlbl0}[c][b]{{\Large Algorithm~\ref{a-sn-chol}: Cholesky (sec.)}}
\psfrag{ylbl0}[c][t]{{\Large Algorithm~\ref{a-sn-projinv}: projected inverse (sec.)}}
\resizebox{0.47\linewidth}{!}{\includegraphics[width=5in]{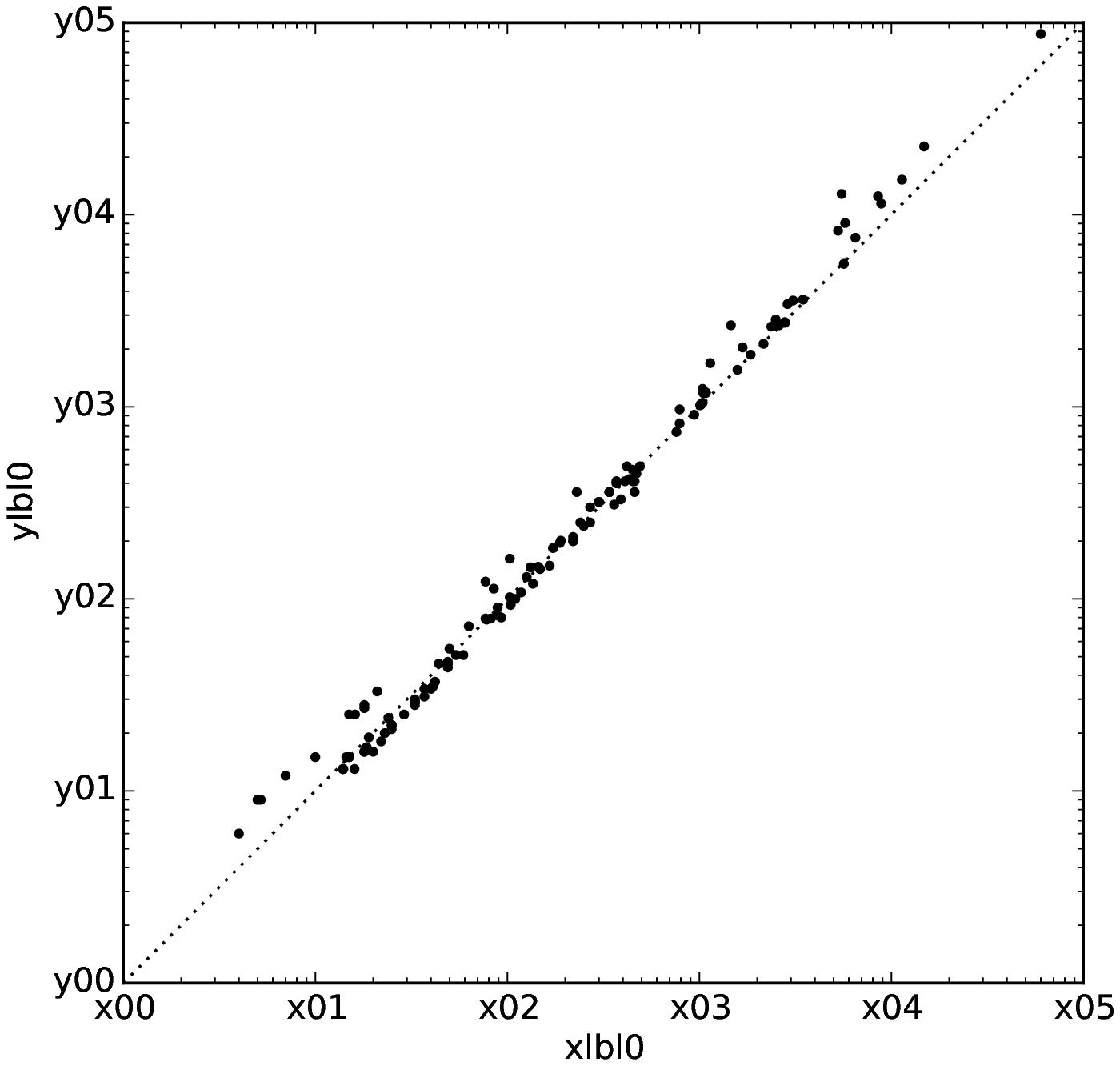}}}
  \ 
  \subfloat[]{\label{fig-ufsmc_mf_parinv-vs-sn_parinv}
\psfrag{x00}[c][c]{${10^{-3}}$}
\psfrag{x01}[c][c]{${10^{-2}}$}
\psfrag{x02}[c][c]{${10^{-1}}$}
\psfrag{x03}[c][c]{${10^{0}}$}
\psfrag{x04}[c][c]{${10^{1}}$}
\psfrag{x05}[c][c]{${10^{2}}$}
\psfrag{x06}[c][c]{${10^{3}}$}
\psfrag{y00}[r][r]{${10^{-3}}$}
\psfrag{y01}[r][r]{${10^{-2}}$}
\psfrag{y02}[r][r]{${10^{-1}}$}
\psfrag{y03}[r][r]{${10^{0}}$}
\psfrag{y04}[r][r]{${10^{1}}$}
\psfrag{y05}[r][r]{${10^{2}}$}
\psfrag{y06}[r][r]{${10^{3}}$}
\psfrag{xlbl0}[c][b]{{\Large Algorithm~\ref{a-projinv}: projected inverse (sec.)}}
\psfrag{ylbl0}[c][t]{{\Large Algorithm~\ref{a-sn-projinv}: projected inverse (sec.)}}
\resizebox{0.47\linewidth}{!}{\includegraphics[width=5in]{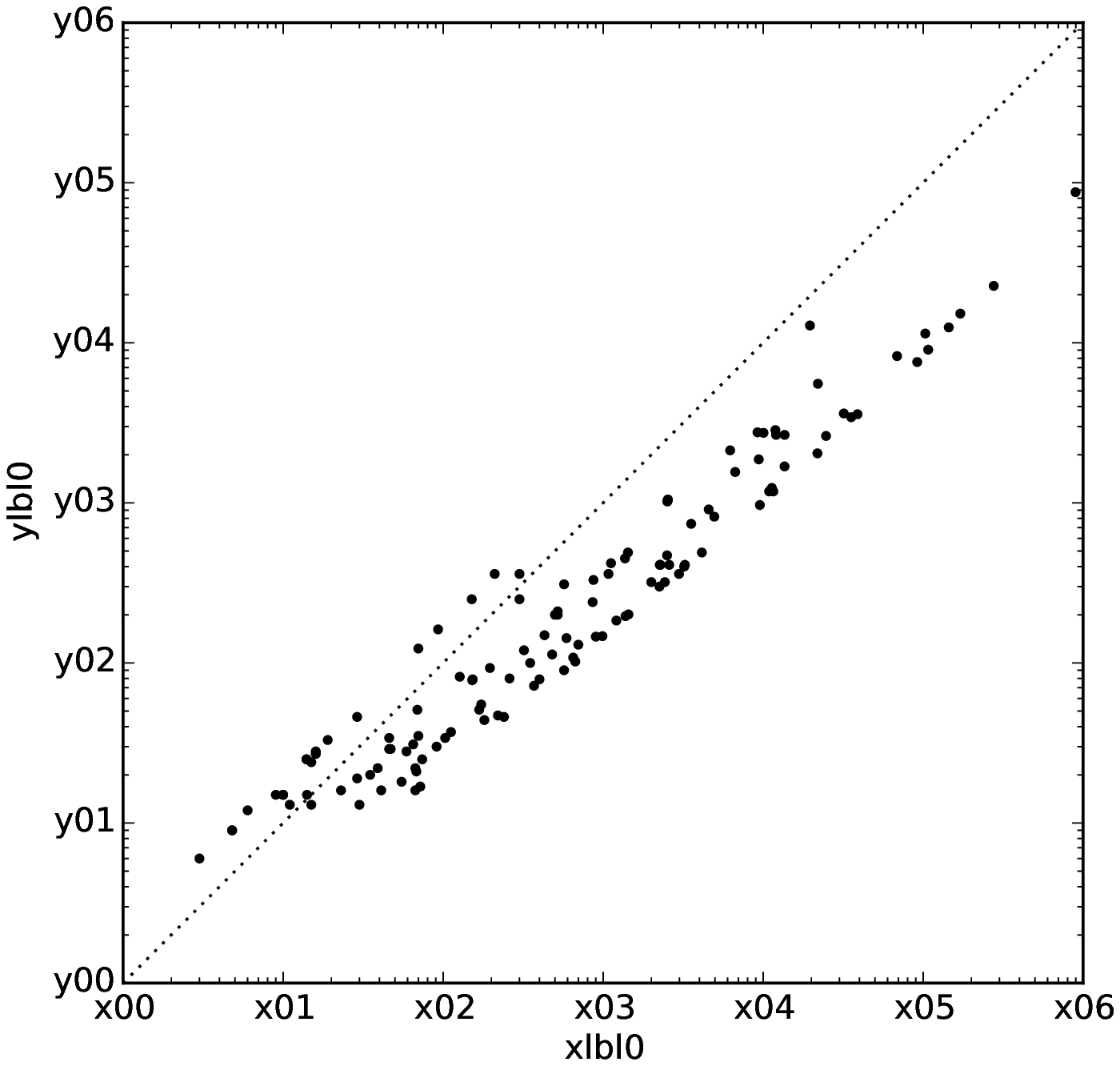}}}
  \caption{Scatter plots of CPU times for supernodal algorithms applied
    to 128 test problems. 
    The plot on the left shows that the cost of computing the Cholesky
    factorization and the projected inverse is approximately the
    same. 
    The plot on the right shows that the supernodal multifrontal
    implementation of the projected inverse algorithm is faster than the
    multifrontal implementation for all but some small problems.}
  \label{fig-supernodal-ufsmc}
\end{figure}

We apply the supernodal multifrontal algorithms 
to the test problems described in Section~\ref{s-numerical-results}. 
The left-hand plot in Figure~\ref{fig-supernodal-ufsmc} shows 
the CPU times for Algorithms~\ref{a-sn-chol} (Cholesky factorizatoin) and 
\ref{a-sn-projinv} (projected inverse or primal gradient). 
The right-hand plot shows shows the CPU times for
Algorithms~\ref{a-sn-projinv} and \ref{a-projinv} (supernodal
multifrontal and multifrontal projected inverse, respectively). 
From the first plot we see that the cost of evaluating the primal
gradient is comparable to the cost of computing the Cholesky
factorization. 
The second plot shows that the supernodal implementation of the primal
gradient is substantially faster than the non-supernodal implementation,
for all but a few of the small problems. 

\section{Conclusions} \label{s-conclusions}
We have derived recursive algorithms for evaluating 
the values, gradients, and Hessians of the primal and dual barriers 
\[
 f(X) = -\log\det X, \qquad
 f_*(S) = \sup_{X\in\symm^n_V} \left(-\Tr(SX) - f(X)\right),
\]
defined for sparse symmetric matrices $X, S\in\SV$ with a given
sparsity pattern $V$, where $V$ is a filled (or chordal) pattern.
Our interest in these algorithms is motivated by their importance 
in interior-point methods for conic optimization with sparse
matrix cone constraints \cite{ADV:10}.  
Similar algorithms can be formulated for closely  
related problems that arise in sparse semidefinite 
programming, for example, the matrix completion techniques used in 
primal-dual methods \cite{FKMN:00,NFFKM:03}.

Our goal was to formulate efficient barrier algorithms based on 
Cholesky factorization techniques for large sparse matrices and, 
specifically, the multifrontal algorithm that has been 
extensively studied in the sparse matrix literature since the 1980s.
The algorithms inherit many of the properties of the multifrontal method. 
This means that a wide range of known techniques from the sparse matrix 
literature can
be used to further improve the algorithms. For example, tree
parallelism and node parallelism are readily exploited in a
multifrontal method \cite{ADL:00}. 
Relaxed super\-nodes and supernode amalgamation techniques 
\cite{DuR:83,AsG:89} have also been
shown to improve the performance. Other improvements include tree
modifications \cite{Liu:88} and memory optimization techniques
\cite{Liu:86b, GuL:06}.

The starting point in this paper was the multifrontal Cholesky
factorization algorithm.
Similar algorithms can be derived from the other popular types 
of sparse factorization algorithms, 
such as the up-looking Cholesky factorization (used in CHOLMOD
\cite{CDHR:08}) or the left-looking Cholesky factorization (a blocked
version of which is used in CHOLMOD's supernodal solver).  It
would be of interest to compare the performance of these algorithms
with the multifrontal algorithms formulated in this paper.

\section*{Acknowledgment}
This material is based upon work supported by the National Science 
Foundation under Grants No.\ ECCS-0824003 and DMS-1115963.
Any opinions, findings, and conclusions or recommendations expressed in  
this material are those of the authors and do not necessarily reflect 
the views of the National Science Foundation.

\newcommand{\etalchar}[1]{$^{#1}$}

\end{document}